\documentclass[12pt]{elsarticle}
\usepackage{amsmath,amsxtra,amssymb,latexsym, amscd,amsfonts,enumerate,ifthen,indentfirst,amstext}
\usepackage{mathrsfs}
\usepackage{longtable}
\usepackage{tablestyles,tabularx, xcolor,multirow}
\usepackage{xcolor}

\usepackage{pgf,tikz,pgfplots}
\pgfplotsset{compat=1.15}
\usepackage{mathrsfs}
\usetikzlibrary{arrows}

\usepackage{smartdiagram}
\usepackage{algorithm2e}

\usepackage[english]{babel}
\usepackage{graphicx}
\usepackage{xcolor,colortbl,color}
\usepackage{graphics,epsfig}
\usepackage{amsthm}
\usepackage{lineno}

\theoremstyle{plain}
\newtheorem{theorem}{\indent\sc Theorem}[section]
\newtheorem{lemma}[theorem]{\indent\sc Lemma}
\newtheorem{corollary}[theorem]{\indent\sc Corollary}
\newtheorem{proposition}[theorem]{\indent\sc Proposition}

\theoremstyle{definition} 
\newtheorem{definition}[theorem]{\indent\sc Definition}
\newtheorem{remark}[theorem]{\indent\sc Remark}

\numberwithin{equation}{section}
\makeatletter
\def\ps@pprintTitle{%
	\let\@oddhead\@empty
	\let\@evenhead\@empty
	\let\@oddfoot\@empty
	\let\@evenfoot\@oddfoot
}
\makeatother

\begin{document}
	
	\begin{frontmatter}
		\title{ Admissible solutions to augmented 
			nonsymmetric $k$-Hessian type equations II. 
			A priori estimates and the Dirichlet problem }  
		
		\author[first]{Tran Van Bang}
		\address[first]{Department of Mathematics,
			Hanoi Pedagogical University 2, 
			 Vinh Phuc, Vietnam }
		\author[second]{Ha Tien Ngoan}
	\address[second]{Institute of  Mathematics, Vietnam Academy of Science and Technology, Hanoi, Vietnam }
		\author[third]{Nguyen Huu Tho}
	\address[third]{Department of  Mathematics, Thuyloi University, Hanoi, Vietnam}
		\author[fourth]{Phan Trong Tien}
	\address[fourth]{Department of Mathematics, Quang Binh University,  Quangbinh, Vietnam}
\begin{abstract}
 Using the established $d$-concavity of the $k$-Hessian type functions $F_k(R)=\log(S_k(R)),$ whose variables are nonsymmetric matrices, we prove $ C^{2, \alpha}(\overline{\Omega}) $  estimates for strictly 
 $(\delta, \widetilde{\gamma}_k) $-admissible solutions to the Dirichlet  
 problem without the well-known regularity condition. A necessary condition for the existence of strictly $\delta$-admissible solutions to the equations is given. By the method of continuity, we provide some sufficient conditions for the unique solvability in the class of strictly  $(\delta,\widetilde{\gamma}_k)$-admissible solutions to the Dirichlet problem, provided that those skew-symmetric matrices in the equations are sufficiently small in some sense.
 
\end{abstract}

	\end{frontmatter} 
	

\section{Introduction}
 This paper is a  continuation of our previous one \cite{2}. We consider the Dirichlet problem for the following nonsymmetric augmented $k$-Hessian type equations
 \begin{equation}\label{1.1}
 		S_{k}\left[D^{2} u-A(x, u, D u)-B(x, u, D u)\right]=f(x, u, D u) \
 	\text { in } \Omega \subset \mathbb{R}^{n}
 \end{equation}
\begin{equation}\label{1.2}
	u(x)=\varphi(x) \ \text { on } \partial \Omega,
\end{equation}
 where $2 \leq k \leq n,$  $\Omega$ is a bounded domain in $\mathbb{R}^{n}$ with smooth boundary $\partial \Omega,$ $D u$ and
 $D^{2} u$ are respectively gradient vector and the Hessian matrix of the unknown function 
 $u: \overline{\Omega} \rightarrow \mathbb{R}, A(x, z, p)=\left[A_{ij}(x, z, p)\right]_{n \times n},$ 
 $B(x, z, p)=\left[B_{ij}(x, z, p)\right]_{n \times n}$ and $f(x, z, p)$ are
 respectively smooth symmetric, skew-symmetric matrices and scalar valued functions, defined on $\mathcal{D}=\overline{\Omega} \times \mathbb{R} \times \mathbb{R}^{n}, \varphi(x)$ is given smooth scalar valued defined on smooth $\partial \Omega.$ We use $x, z, p, R$ to denote points in $\overline{\Omega}, \mathbb{R}, \mathbb{R}^{n}$ and $\mathbb{R}^{n \times n}$
 respectively. Here,
 \begin{equation*}
 	 S_{k}(R)=\sigma_{k}(\lambda(R)),
 \end{equation*}
  where $\lambda(R)=\left(\lambda_{1}, \lambda_{2}, \cdots, \lambda_{n}\right) \in \mathbb{C}^{n}$ is the vector
 of eigenvalues of the matrix $R=\left[R_{ij}\right]_{n \times n} \in \mathbb{R}^{n \times n},$
\begin{equation*}
	\sigma_{k}(\lambda)=\sum_{1\le i_1<\cdots <i_k\le n} \lambda_{i_1} \lambda_{i_{2}} \cdots \lambda_{i_{k}}
\end{equation*}
is the elementary symmetric polynomial of degree $k.$ Noting that, since $R \in \mathbb{R}^{n \times n}, S_{k}(R)$ is real-valued. When $B(x, z, p) \equiv 0$ the equations $\eqref{1.1}$ become
\begin{equation}\label{1.5}
	S_{k}(D^{2} u-A(x, u, D u))=f(x, u, D u) \text{ in } \Omega,
\end{equation}
which are symmetric augmented $k$-Hessian type equations. When $k=n$ the equations \eqref{1.5} are the Monge-Ampère type equations:
\begin{equation*}
	\det\left(D^{2} u-A(x, u, D u)\right)=f(x, u, D u)  \text{ in } \Omega,
\end{equation*}
the Dirichlet problem for which had been studied in \cite{3}-\cite{6}, \cite{8}, \cite{13}-\cite{15}, \cite{17}.

For $u(x) \in C^{2}(\overline{\Omega})$ and $x \in \overline{\Omega}$ we set 
\begin{equation}\label{1.7} \omega(x, u) =D^{2} u(x)-A(x, u(x), D u(x))    =\left[\omega_{i j}(x, u)\right]_{n \times n} .
 \end{equation}
For  $1 \leq k \leq n$  we denote by $ \Gamma_{k}$ the following cone in $\mathbb{R}^{n}:$
\begin{equation*}
	\Gamma_{k}=\left\{\lambda \in \mathbb{R}^{n}: \sigma_{j}(\lambda)>0, j=1, \cdots, k\right\}.
\end{equation*}
 When $k=n$  we have 
 $$\Gamma_{n }=\left\{ \lambda=(\lambda_{1},\cdots, \lambda_{n})\in \mathbb{R}^n; \lambda_{j}>0,j=1,\cdots,n \right\}.$$
 
The equations \eqref{1.5} have been considered in \cite{7}, \cite{9}, \cite{16}, \cite{18}. A function $u(x) \in C^{2}(\overline{\Omega})$ is said to be an admissible solution (\cite{9}) to the equation \eqref{1.5} if $\lambda(\omega(x, u)) \in \Gamma_{k}$ for any $x \in \overline{\Omega}.$ Under the assumption of regularity condition (see \eqref{T7.16}) for the matrix $A(x, z, p)$ and that of existence of an admissible subsolution, by using the concavity of the function $\sqrt[k]{\sigma_{k}(\omega)},$ the authors of \cite{9} has proved the unique existence of admissible solution to the Dirichlet problem for \eqref{1.5}.

 The nonsymmetric Monge-Ampère type equations
\begin{equation}\label{1.10}
\det\left(D^{2} u-A(x, u, D u)-B(x, u, D u)\right)=f(x, u, D u)	
\end{equation}
has been considered in \cite{11}, \cite{12}. The main difficulty in this case is that the both functions $\sqrt[n]{\det R}$ and $\log \left(\det R \right)$ are not concave. To overcome this difficulty, the following class of elliptic solutions to \eqref{1.10} are introduced as follows.
\begin{definition}[\cite{11}, \cite{12}]
Suppose $u(x) \in C^{2}(\overline{\Omega}).$ Then
\begin{itemize}
	\item[(i)] The function $u(x)$ is said to be an elliptic solution to \eqref{1.10} if the following condition holds
	\begin{equation}\label{1.11}
		\lambda_{u}:=\inf\limits_{x\in \overline{\Omega}} \lambda_{\min }(\omega(x, u))>0,
	\end{equation}
where $\lambda_{\min }(\omega)$ is the least eigenvalue of $\omega;$
	\item[(ii)]  Suppose $0<\delta<1.$ The function $u(x)$ is said to be $\delta$-elliptic solution to \eqref{1.10} if it is elliptic one and it holds

	\begin{equation}\label{1.12}
	\mu(B) \leq \delta \lambda_{u},
	\end{equation}

where the matrix $B(x, z, p)$ is assumed to belong to $B C(\mathcal{D})$ and
\begin{equation*}
	\mu(B):=\sup _{\mathcal{D}}\|B(x, z, p)\|,
\end{equation*} 
here $\|B\|$ stands for the operator norm of the matrix $B.$
\end{itemize}
\end{definition}

For $u(x)\in C^2(\overline{\Omega})$ and $x\in\overline{\Omega}$ we set
\begin{equation}\label{1.7a}
R(x,u)=\omega(x,u)-B(x,u,Du)=[R_{i j}(x,u)]_{n\times n},
\end{equation}
where $\omega(x,u)$ is defined by \eqref{1.7}. 

In connection with the $\delta$-elliptic solutions, the following convex and unbounded set of nonsymmetric matrices $R$ had been introduced for   $0<\delta<1, \mu>0$  (\cite{11})
 \begin{equation*}
 	D_{\delta, \mu}=\left\{R \in \mathbb{R}^{n \times n}: R=\omega+\beta, \omega^{T}=\omega,
 \beta^{T}=-\beta, \omega>0,\|\beta\| \leq\mu, \mu \leq \delta \lambda_{\min }(\omega)\right\} 
 \end{equation*}
as a domain for $F(R)=\log (\det R).$ We note that if $u(x)$ is a  $\delta$-elliptic solution then $R(x,u)\in D_{\delta, \mu(B)}$ for any $x\in\overline{\Omega }.$ The notion of $d$-concavity for the function $F(R)=\log (\det R)$ for $d \geq 0$ had been introduced in \cite{11} as follows:

\begin{definition} The function $F(R)$ is said to be $d$-concave on $D_{\delta, \mu}$ if for any $R^{(0)}=\left[R^{(0)}_{ij}\right]=\omega^{(0)}+\beta^{(0)}, R^{(1)}=\left[R_{i j}^{(1)}\right]=\omega^{(1)}+\beta^{(1)} \in D_{\delta,\mu}$ the following inequality holds: 
\begin{equation*}
	F\left(R^{(1)}\right)-F\left(R^{(0)}\right) \leq\sum_{i, j=1}^{n} \frac{\partial F\left(R^{(0)}\right)}{\partial R_{ij}}\left(R_{ij}^{(1)}-R_{i j}^{(0)}\right)+C.\frac{\left| \beta^{(1)}-\beta^{(0)} \right|^2}{\lambda_{\min }^2\left(\omega^{(\tau)}\right)},
\end{equation*}	
where $\omega^{(\tau)}=(1-\tau)\omega^{(0)}+\tau\omega^{(1)}, 0<\tau<1.$ 
\end{definition}

The $d$-concavity of the function $F(R)=\log (\det R)$ had been established (\cite{11}, Theorems 2 and 3), where $C$ depends only on $\delta, n$ and does not depend on $\mu.$ Then the $d$-concavity, the regularity condition for the matrix $A(x,z,p)$ and the assumption on existence of an elliptic subsolution $\underline{u}(x)$ to the problem \eqref{1.10}-\eqref{1.2} with $B(x,z,p)=0,$ enable to get 
	$C^{2, \alpha}(\overline{\Omega})$-estimates for $\delta$-elliptic solutions to the Dirichlet problem \eqref{1.10}-\eqref{1.2} with some $0<\alpha<1$ and then to get the solvability of the problem (\cite{12}, Theorems 3 and 4). 
	
	In this paper, for  the cases $2 \leq  k \leq  n$ we prefer to replace the notions of elliptic and $\delta$-elliptic solutions respectively by 
	the notions of strictly admissible and strictly $\delta$-admissible solutions to the equations \eqref{1.1} that are defined respectively as the same as elliptic and $\delta$-elliptic solutions for the Monge-Ampère type equations \eqref{1.10}. But to get the $d$-concavity of the functions $F_k(R)=\log(S_k(R))$ we have to restrict more on these classes of strictly admissible solutions. To do this we define a subcone $\sum_{(\widetilde{\gamma}_k)}$ in $\Gamma_{n}, 0<\tilde{\gamma}_k<1,$ as follows.
	
\begin{definition}[\cite{2}]\label{dn1.3}
Suppose $1\le k\le n.$ The subcone $\Gamma_{(\widetilde{\gamma}_k)}$ consists of all 
$\lambda=\left(\lambda_{1},\cdots, \lambda_{n}\right) \in \Gamma_{n} ,$ such that
\begin{equation*}
	\lambda_{\min } \geq \widetilde{\gamma}_{k} \lambda_{\max },
\end{equation*}
where  $\lambda_{\min }=\min_{1\le j\le n} \lambda_{j}, $ $\lambda_{\max }=\max _{1\le j\le n} \lambda_{j} $ and $\widetilde{\gamma}_{k} $ is chosen appropriately in each concrete
 problem and satisfies the following conditions:
\begin{itemize}
	\item[(i)]  If $k\in \{2,3,n-1,n\},$ then $0<{\gamma}_{k}< \widetilde{\gamma}_{k} <1,$ where $\gamma_{k}$ is a some positive number that is less than 1 and must be also determined in each case;
	\item[(ii)]   If $\left[\frac{n}{2}\right]+1 \leq k \leq n-2,$ then 
	\begin{equation}\label{1.16}
		\gamma_{k}=\frac{n-k}{k}<\widetilde{\gamma}_{k} <1;
	\end{equation}
	\item[(iii)]   If $4 \leq  k \leq \left[\frac{n}{2}\right], $ then 
	\begin{equation}\label{1.17}
		\gamma_k=\gamma_{n-k+2}=\frac{k-2}{n-(k-2)}< \widetilde{\gamma}_{k}<1,
	\end{equation}
\end{itemize}
where $\gamma_k,$ $2 \leq  k \leq n-1,$ have been already defined in \cite{2} as above.
\end{definition}
Now the domain of the function
$F_{k}(R)=\log \left(S_{k}(R)\right)$ is introduced as follows:
\begin{definition}[\cite{2}] Suppose $0<\delta<1, \mu>0$
and $0<{\gamma}_{k}< \widetilde{\gamma}_{k}<1$ that have been defined as above. We set
\begin{equation*}
D_{\delta, \mu, \widetilde{\gamma}_{k}}=\left\{R=\omega+\beta \in D_{\delta, \mu }; \lambda(\omega) \in \Sigma_{\left(\widetilde{\gamma}_{k}\right)}\right\}.
\end{equation*}
Noting that all the sets $\Sigma_{\left(\widetilde{\gamma}_{k}\right)},$ $D_{\delta, \mu }$ and $D_{\delta, \mu, \widetilde{\gamma}_{k}}$
are convex and unbounded.
\end{definition}
We recall now some following results from \cite{2} (Theorem 1) and \cite{11} (Proposition 5.1 and Theorem 1.6)
for the functions $F_{k}(R)=\log \left(S_{k}(R)\right).$
\begin{theorem}[\cite{2}, \cite{11}]\label{dl1.5}
	 Suppose $2 \leq  k \leq  n$ and $0<\gamma_{k}< \widetilde{\gamma}_k<1$ are defined as in Definition \ref{dn1.3}. Then there exist $\delta_{k}, 0<\delta_{k}<1, \delta_{k}=\delta_{k}\left(k, n, \widetilde{\gamma}_{k}\right)$ 
	 if $2\le k\le (n-1)$ and $\delta_k$ may be any positive number that is less than 1 when $k=n$
	 and $C_{j}>0,$
	$C_j=C_j\left(k, n, \widetilde{\gamma}_{k},\delta_k\right)$ such that for all $\delta, 0<\delta <  \delta_{k}$ and 
	\begin{itemize}
		\item [(i)] for all $R=\omega+\beta\in  D_{\delta, \mu, \widetilde{\gamma}_{k}},$  $M=P+Q \in \mathbb{R}^{n \times n}, P^{T}=P, Q^{T}=-Q$ the following estimates hold
		\begin{equation}\label{1.19}
			d^{2} F_{k}(R, P) \leq -C_{1} \frac{|P|^{2}}{\lambda_{\max }^{2}(\omega)},
		\end{equation}
	\begin{equation}\label{1.20}
		d^{2} F_{k}(R, M) \leq  C_{2} \frac{|Q|^{2}}{\lambda_{\min }^{2}(\omega)},
	\end{equation}
where for $M=[M_{ij}]\in \mathbb{R}^{n\times n}	,$ $|M|^2=\sum_{i,j=1}^{n}|M_{ij}|^2;$ 
\item[(ii)] for all $R^{(0)}, R^{(1)}\in D_{\delta, \mu, \widetilde{\gamma}_{k}},$ $R^{(0)}=\omega^{(0)}+\beta^{(0)},$ $R^{(1)}=\omega^{(1)}+\beta^{(1)},$ the following $d$-concavity of the function $F_k(R)$, that is a consequence of \eqref{1.20}, holds
\begin{equation}\label{1.21}
F_{k}\left(R^{(1)}\right)-F_{k}\left(R^{(0)}\right)  \leq\sum_{i, j=1}^{n} \frac{ \partial F_{k}\left(R^{(0)}\right)}{\partial R_{ij}}\left(R_{ij}^{(1)}-R_{i j}^{(0)}\right) 
		+C_{2} \frac{|\beta^{(1)}-\beta^{(0)}|^{2}}{\lambda_{\min }^{2}\left(\omega^{(\tau)}\right)}, 
\end{equation}	
where $\omega^{(\tau)}=(1-\tau) \omega^{(0)}+\tau \omega^{(1)}, 0<\tau<1.$
		\end{itemize}
\end{theorem}

From here and throughout the paper we always assume  that the parameters $\delta, \widetilde{\gamma}_k$ are defined as follows:
\begin{equation}\label{1.12a}
	0<\gamma_{k} <\widetilde{\gamma}_k <1, \quad 0<\delta <\delta_k <1,
\end{equation}
where $0<\gamma_{k} <1$ have been defined in Definition \ref{dn1.3} and $0<\delta_k<1$ has been determined in  Theorem \ref{dl1.5}.

	Now other types of strictly admissible solutions to the equations \eqref{1.1} are introduced.
 \begin{definition}
 Suppose $u(x) \in C^{2}(\overline{\Omega})$ is a
 	strictly admissible solution to \eqref{1.1}, $0<\delta<1, 0< \widetilde{\gamma}_{k}< 1$ as above, then
 	\begin{itemize}
 	\item [(i)]  It is said to be strictly $\widetilde{\gamma}_{k}$-admissible solution to \eqref{1.1} if 
 	\begin{equation}\label{1.22}
 		\gamma_u:=\inf_{x\in\overline{\Omega} }\left[ \frac{\lambda_{\min}(\omega(x,u))}{\lambda_{\max}(\omega(x,u))}\right] \ge \widetilde{\gamma}_{k},
 	\end{equation}
 	where $\omega(x,u)$ is defined by \eqref{1.7}, i.e. $\lambda(\omega(x,u))\in \sum_{(\widetilde{\gamma}_k)}, \forall x\in \overline{\Omega}; $
 	\item[(ii)]  It is said to be strictly $(\delta, \widetilde{\gamma}_{k})$-admissible solution to \eqref{1.1} if it is both strictly $\delta$-admissible and strictly $\widetilde{\gamma}_{k}$-admissible solution to \eqref{1.1}, i.e. \eqref{1.12}, \eqref{1.22} hold and $R(x,u)\in D_{\delta, \mu(B), \widetilde{\gamma}_k}, \forall x\in \overline{\Omega}.$
 \end{itemize}
 	 	 \end{definition}
  	 \begin{remark}
  	 	 The condition \eqref{1.22} seems to be rather strict one, because the equation \eqref{1.1} becomes indeed uniformly elliptic at solutions of this kind and the $C^2(\overline{\Omega })$-estimates for solutions are easily obtained.       But the condition \eqref{1.22} is actually needed, because it allows the function $F_k(R(x,u))=\log \left(S_k(R(x,u))\right)$ to be $d$-concave with respect to $R(x,u)$ for $x\in \overline{\Omega },$ with the aid of which one can prove the Holder continuity of $D^2 u(x)$ in $\overline{\Omega}$. The condition \eqref{1.22} is only a structural one for solutions of the problem \eqref{1.1}-\eqref{1.2}, but it is not structural one for the equations \eqref{1.1}. The most important structural conditions for the data  $A(x,z,p), f(x,z,p)$ and $\Omega$, as it will be clear later in an example at the last section of the paper, must be those ones, under which there exists a strictly $\widetilde{\gamma}_k$-admissible subsolution $\underline{u}(x)$ of the problem \eqref{1.1}-\eqref{1.2}. 
  	 	 
  	 \end{remark}
\begin{remark}
When $k=n,$ in \cite{8} the authors did not assume the uniform ellipticity condition for elliptic solutions to the Monge-Ampère type equations. But the elliptic solution $u(x),$ that exists and is unique in \cite{8}, is actually a strictly $\widetilde{\gamma}_n$-admissible one, where $\widetilde{\gamma}_n$ is some positive number, that is less than $1.$ Indeed, thanks to assumptions on regularity condition \eqref{T7.16} and some additional structural conditions on $A(x,z,p), f(x,z,p)$ and the assumption on existence of elliptic subsolution $\underline{u}(x)$ to the problem	
\eqref{1.1}-\eqref{1.2}, the authors had proved that there exist $M_0>0,M_1>0,M_2>0,$ $0<M_3<M_4$ such that 
\[\underset{x\in\overline{\Omega}}{\sup}|u(x)|\leq M_0,\quad \underset{x\in\overline{\Omega}}{\sup}|Du(x)|\leq M_1, \quad \underset{x\in\overline{\Omega}}{\sup}|D^2u(x)|\leq M_2,\]
from which one obtains 
\[\underset{x\in\overline{\Omega}}{\inf}\lambda_{\min}(\omega(x,u))\geq M_3,\quad \underset{x\in\overline{\Omega}}{\sup}\lambda_{\max}(\omega(x,u))\leq M_4,\]
and therefore \eqref{1.22} follows with $\tilde{\gamma}_n=\frac{M_3}{M_4}.$
\end{remark}

The purpose of the paper is to study the solvability of the problem \eqref{1.1}-\eqref{1.2} in the class of strictly $\left(\delta, \widetilde{\gamma}_{k}\right)$-admissible solutions without regularity condition for the matrix $A(x,z,p).$ The paper is organized as follows. In Section 2 we establish the comparison principle (Theorem \ref{dl2.2}) for strictly $\delta$-admissible solutions. This principle is analogous to that for Monge-Ampère type \eqref{1.10} equations (\cite{11}).
 In Section 3,  for strictly $\left(\delta, \widetilde{\gamma}_{k} \right)$-admissible solutions to the Dirichlet problem \eqref{1.1}-\eqref{1.2}, we estimate eigenvalues of the matrices $\omega(x,u)$ at any $x\in \overline{\Omega}.$ It is interesting that for this kind of solutions, to do this, we do not need neither $d$-concavity of the function $F_k(R),$ nor regularity condition for the matrix $A(x,z,p).$
\begin{proposition}\label{md1.7}
Suppose $A(x, z, p), f(x, z, p) \in C(\mathcal{D}),$ $B(x, z, p) \in B C(\mathcal{D}), f(x, z, p)>0.$ Suppose $u(x) \in C^{2}(\overline{\Omega})$
is a strictly admissible solution to the equation \eqref{1.1} and there exist $M_{0}>0, M_{1}>0$ such that
$$
\sup_{\overline{\Omega}}|u(x)| \leq M_{0}, \quad \sup _{\overline{\Omega}}|D u(x)| \leq M_{1}.$$
We set
\begin{equation}\label{1.22a}
		f_{0}=\inf_{x \in \overline{\Omega} \atop
	|z|\le M_0, |p|\le M_1} f(x, z, p),\quad   f_{1}=\sup_{x \in \overline{\Omega} \atop
	|z|\le M_0, |p|\le M_1} f(x, z, p).
\end{equation}

Then the following assertions hold for any $x \in \overline{\Omega}:$
\begin{itemize}
	\item[(i)] \begin{equation}\label{1.23}
		\quad 0<\lambda_{\min }(\omega(x, u)) \leq\left[\frac{f_{1}}{ \binom{n}{k}} \right]^{\frac{1}{k}};
	\end{equation}
\item[(ii)]  If $u(x)$ is a strictly $\delta$-admissible solution, then besides \eqref{1.23}, the following inequality is true
\begin{equation}\label{1.24}
	\left[\frac{   (1+\delta^{2} )^{-\left[\frac{k}{2}\right]} f_{0} }{\binom{n}{k} }\right]^{\frac{1}{k}} \leq\lambda_{\max }(\omega(x, u));  
\end{equation}

\item[(iii)]  If $u(x)$ is a strictly $ \widetilde{\gamma}_{k}$-admissible solution, then besides \eqref{1.23}, the following inequality is true
\begin{equation}\label{1.25a}
	\lambda_{\max }(\omega(x, u))
	\le \frac{1}{\widetilde{\gamma}_{k}}	\left[\frac{f_1}{\binom{n}{k}}\right]^{\frac{1}{k}};
\end{equation}
\item[(iv)] If $u(x)$ is a strictly $(\delta,\widetilde{\gamma}_k)$-admissible solution, then
\begin{equation}\label{1.25}
	\widetilde{\gamma}_{k} \left[\frac{   (1+\delta^{2} )^{-\left[\frac{k}{2}\right]} f_{0} }{\binom{n}{k}}\right]^{\frac{1}{k}} \leq\lambda_{\min }(\omega(x, u))
	\le 	\left[\frac{f_1}{\binom{n}{k}}\right]^{\frac{1}{k}};
\end{equation}
\begin{equation}\label{1.26}
	\left[\frac{   (1+\delta^{2} )^{-\left[\frac{k}{2}\right]} f_{0} }{\binom{n}{k}}\right]^{\frac{1}{k}} \leq\lambda_{\max }(\omega(x, u))
	\le \frac{1}{\widetilde{\gamma}_{k}}	\left[\frac{f_1}{\binom{n}{k}}\right]^{\frac{1}{k}}.
\end{equation}
\end{itemize}
\end{proposition}

Under some structure conditions on the matrix $A(x,z,p),$ proposed by N.S Trudinger and his colleagues in \cite{8}, by using the comparison principle (Theorem \ref{dl2.2}) we obtain $C^2(\overline{\Omega})$-estimates for strictly $(\delta,\tilde{\gamma}_k)$-admissible solutions $u(x)$ in the following theorem.
\begin{theorem}\label{dl1.8}
	Assume that $0< \delta<1,$ $0< \widetilde{\gamma}_k <1$ are defined as in  \eqref{1.12a} and the following conditions fulfill:
	\begin{itemize}
 \item[(i)] $A(x, z, p)\in C^3(\mathcal{D})$ and satisfies structure conditions:
  $$A(x, z, p) \geq-\gamma_{0}\left(1+|p|^{2}\right) E_{n},\  \gamma_{0}>0,$$
 $$\lambda_{\max }(A(x, z, 0)) \geq 0, \quad D_{z} A(x, z, p) \geq 0;$$
 \item[(ii)] $f(x, z, p)\in C^3(\mathcal{D})$ and
 $$f(x,z,p)>0 \text{ in } \mathcal{D},$$
 $$\inf_{\mathcal{D}}\left[\frac{D_z f(x,z,p)}{f(x,z,p)} \right]\ge \frac{k\delta}{(1+\delta^2)}\beta_1, \ \beta_1>0,$$
  \item[(iii)] There exists a strictly $\widetilde{\gamma}_{k}$-admissible subsolution $\underline{u}(x)$ to the problem
  \begin{align}
 	&S_{k}\left(D^{2} u-A(x, u, D u)\right) =f(x, u, D u) \  \text { in } \Omega, \label{1.27} \\
 	&u=\varphi \ \text { on } \partial \Omega; \label{1.28} 
 \end{align}
\item[(iv)]  Suppose  $u(x) \in C^{2}(\overline{\Omega})  $ is a strictly $\left(\delta,\widetilde{\gamma}_k\right) $-admissible solution to the problem \eqref{1.1}-\eqref{1.2};
\item[(v)] $B(x, z, p)\in BC^3(\mathcal{D})$ and
$$ \mu(B) \leq \delta \min \left(\lambda_{u}, \lambda_{\underline{u}}\right),$$
$$ \mu\left(D_{z} B\right) \leq \beta_{1} \min\left(\lambda_{u}, \lambda_{\underline{u}}\right). $$
Then there exist  $M_{0}>0, M_{1}>0, C_{3}>0,$ that depend only on $ \delta, k, n, \widetilde{\gamma}_{k}, \beta_{1}, A, \underline{u}, f, \varphi$   
such that
$$\sup_{\overline{\Omega}}|u(x)|\le M_0,\quad \sup_{\overline{\Omega}}|Du(x)|\le M_1$$
and
\begin{equation}\label{1.29}
	\|u\|_{C^2(\overline{\Omega})}\le C_3.
\end{equation}
	\end{itemize}
\end{theorem}
Using \eqref{1.29}, the ellipticity estimate \eqref{3.6} and the $d$-concavity in the sense of  \eqref{1.21} of the function $F_k(R)=\log \left(S_k(R)\right),$ at the end of Section 3 we show the H\"older continuity of second-order derivatives $D^2u$ with some $0<\alpha <1$ inside $\Omega,$ if $u(x)\in C^4(\Omega).$ 
Here, besides the quantities, on which $C_3$ depends, $\alpha$ depends also on $\mu_2(B),$ where
\begin{equation}\label{1.23a}
\mu_2(B)=\|B(x,z,p)\|_{BC^2(\mathcal{D})}.
\end{equation}

 In Section 4 we consider $x^0\in \partial\Omega.$
By translation and rotation, we can assume that $x^{0}$ is the origin of coordinates and the unit inner normal at $x^0$ is on the axis $O x_n.$
Suppose that in a neighborhood $\mathcal{N}  $ of $x^{0},$ the boundary $\partial\Omega$ is the graph of the function
$$
x_{n}=h\left(x'\right), \quad x'=\left(x_{1}, \ldots, x_{n-1}\right),
$$
where $h(x')\in C^4$ and
\begin{equation*}
	h\left(0'\right)=0, \quad D h\left(0'\right)=0.
\end{equation*}
We change $x=\left(x_{1}, \ldots, x_{n}\right)$ into $y=\left(y_{1},\ldots, y_{n}\right)$ by the mapping
\begin{equation}\label{1.31}
	y=\psi(x)=\left(\psi_{1}(x), \ldots, \psi_{n}(x)\right)=\left(x', x_{n}-h\left(x'\right)\right),
\end{equation}
where $y\in \widetilde{\mathcal{N} }:=\psi(\mathcal{N} ). $ From \eqref{1.31}
we have
\begin{equation}\label{1.32}
	x=\widetilde{\psi}(y)=\left(\widetilde{\psi}_{1}(y), \ldots, \widetilde{\psi}_{n}(y)\right)=\left(y', y_{n}+h\left(y'\right)\right).
\end{equation}
We set
\begin{equation}\label{1.33}
	v(y)=u(x)-\varphi(x), 
\end{equation}
\begin{equation*}
	J(x)=\frac{D \psi(x)}{D x}=\left[\frac{\partial \psi_{i}(x)}{\partial x_{j}}\right]_{n \times n}= \left[
	\begin{array}{ccc}
		\frac{\partial \psi_{1}(x)}{\partial x_1 } & \cdots & \frac{\partial \psi_{1}(x)}{\partial x_n }   \\
		 \cdots &  \cdots &  \cdots  \\
		\frac{\partial \psi_{n}(x)}{\partial x_1 } & \cdots & \frac{\partial \psi_{n}(x)}{\partial x_n }  \\
	\end{array}
	\right],
\end{equation*}
  where $x$ and $y$ are related by \eqref{1.31}  and \eqref{1.32}. In \eqref{1.33} we assume that the function $\varphi(x)$ has been extended smoothly from $\partial \Omega$ into some neighborhood of $\partial \Omega.$

 Then we have 
 \begin{equation}\label{1.35}
 	J(x)= \left[\begin{array}{cc}
 		E_{(n-1)} &  0 \\
 		-D h(x') & 1
 	\end{array}\right]= J(y),\quad  J^{T}(x)=\left[\begin{array}{cc}
 		E_{(n-1)} & -(D h)^T(x') \\
 		0 & 1
 	\end{array}\right]=J^{T}(y),
 \end{equation}
 \begin{equation}\label{1.27a}  
 	J^{-1}(x) =\left[\begin{array}{cc}
 		E_{(n-1)} & 0 \\
 		D h(x') & 1
 	\end{array}\right]=J^{-1}(y),\  {(J^{-1})}^{T}(x)=\left[\begin{array}{cc}
 		E_{(n-1)} & (D h)^T(x') \\
 		0 & 1
 	\end{array}\right] =\left(J^{-1}\right)^T(y),
 \end{equation}
where $E_{(n-1)}$ is the unit matrix of size $(n-1)$ and $Dh$ stands for the
row vector  
$$Dh(x')=(D_{x_1}h(x'),\cdots,D_{{x_{n-1}}}h(x'))= (D_{y_1}h(y'),\cdots,D_{{y_{n-1}}}h(y'))=Dh(y') .$$

We have
\begin{equation}\label{1.37}
	D_x u= (D_y v)J+D \varphi,
\end{equation}
\begin{equation}\label{1.38}
	D_x^{2} u=J^{T} D_y^{2} v J+\sum_{m=1}^{n} D_{y_m}v D^{2} \psi_{m}+D_x^{2} \varphi,
\end{equation}
where $D u=\left(D_{x_1} u, \cdots, D_{x_n} u\right), D v= (D_{y_1} v, \cdots, D_{y_n} v).$
We set further on the base of \eqref{1.37}, \eqref{1.38}:
\begin{equation}\label{1.39}
	\left\{\begin{array}{l}
		\widetilde{A}(y, z, p)=\left(J^{-1}\right)^T\left[A(\widetilde{\psi}(y), z+\varphi(\widetilde{\psi}(y)), pJ +D_x \varphi(\widetilde{\psi}(y))) 
		\right.\\
	\hspace*{4cm}	\left.-\sum_{m=1}^{n} p_{m} D_x^{2} \psi_{m}(\widetilde{\psi}(y))-D_x^{2} \varphi(\widetilde{\psi}(y))\right] J^{-1}, \\
		\widetilde{B}(y, z, p)=\left(J^{-1}\right)^T B(\widetilde{\psi}(y), z+\varphi(\widetilde{\psi}(y)), pJ +D_x \varphi(\widetilde{\psi}(y)))\left(J^{-1}\right)^{T}, \\
		\widetilde{f}(y, z, p)=f(\widetilde{\psi}(y), z+\varphi(\widetilde{\psi}(y)), pJ+D_x \varphi (\widetilde{\psi}(y))).
	\end{array}\right.
\end{equation}
It follows from \eqref{1.38}-\eqref{1.39} that 
\begin{equation*}
	D^{2} u-A(x, u, D u)-B(x, u, Du)= J^{T}\left[D^{2} v-\widetilde{A}(y, v, D v)-\widetilde{B}(y, v, D v)\right] J.
	\end{equation*}
The equation \eqref{1.1} becomes 
\begin{equation}\label{1.41}
	S_{k}\left[J^{T}\left(D^{2} v-\widetilde{A}(y, v, D v)-\widetilde{B}(y, v, D v)\right) J\right]=\widetilde{f}\left(y, v, D v\right) \text { in } \widetilde{\Omega}_{\rho},
\end{equation}
where $\widetilde{\Omega}_{\rho}=\left\{\left(y', y_{n}\right):|y|<\rho, y_{n}>0\right\},$ $v(y)$ satisfies condition:
\begin{equation}\label{1.31a}
	v\left(y', y_{n}\right)=0  \text { when } y_{n}=0,\left|y'\right|<\rho,\ \rho>0.
\end{equation}
We set for $v(y)\in C^2(\overline{\widetilde{\Omega}_{\rho}})$
\begin{equation*}
	\widetilde{\omega}(y,v)=D^2v-\widetilde{A}(y, v, D v),
\end{equation*}
\begin{equation*}
	\widetilde{R}(y,v)=\widetilde{\omega}(y,v)- \widetilde{B}(y, v, D v)=\left[\widetilde{R}_{ij}\right]_{n\times n}.
\end{equation*}
Suppose $i_{1} i_{2} \cdots i_{k}$ and $j_{1} j_{2} \cdots j_{k}$ are indices such that
$$
1 \leq i_{1}<i_{2}<\cdots<i_{k} \leq n, \quad 1 \leq j_{1}<j_{2}<\cdots<j_{k} \leq n.
$$
We denote
\begin{equation*}
	\widetilde{R}^{(k)}_{i_1\cdots i_k,j_1\cdots j_k}=\left[ \widetilde{R}_{i_p j_q}\right]^k_{p,q=1}.
\end{equation*}

\begin{proposition}\label{md1.9} In a neighborhood of the origin $y^{0}=0$ the equation \eqref{1.41} can be rewritten in  the form:
\begin{equation}\label{1.46}
	S_{k}(\widetilde{R}(y, v))+H_{k}\left(y', \widetilde{R}(y, v)\right)=\widetilde{f}(y, v, D v),\  y\in \widetilde{\Omega_{\rho}},
\end{equation}
where $H_n(y',\widetilde{R})=0$ and if $2\leq k\leq n-1$ then
 \begin{align}\label{1.47}
 	\begin{split}
 &H_k(y',\widetilde{R}) =\sum_{1\le i_1<\cdots <i_{k-1}\le n-1} \left[ \sum_{m\in\{1,\cdots,n-1\}\backslash \{i_1,\cdots,i_{k-1}\}} (D_m h(y'))^2 \right] \det \widetilde{R}^{(k)}_{i_1\cdots i_{k-1} n, i_1\cdots i_{k-1} n}\\
	&+ (-1)^{k-1}\sum_{1\le i_1<\cdots <i_{k}\le n-1\atop 1\le j_1<\cdots <j_{k-1}\le n=j_k}\left[\sum_{m=1}^{k}(-1)^m (D_{i_m} h(y')) \delta_{i_1j_1}\cdots \delta_{i_{m-1}j_{m-1}}  \delta_{i_{m+1}j_{m}} \cdots \delta_{i_{k}j_{k-1}}\right]  \\
	&\hspace*{3cm} \times \left(\det \widetilde{R}^{(k)}_{i_1\cdots i_k,j_1\cdots j_{k-1}n}+\det   \widetilde{R}^{(k)}_{j_1\cdots j_{k-1}n,i_1\cdots i_{k}}   \right)\\
	& - \sum_{1\le i_1<\cdots <i_{k-1}\le n-1\atop {1\le j_1<\cdots <j_{k-1}\le n-1 \atop \{i_1,\cdots,i_{k-1}  \} \ne \{j_1,\cdots,j_{k-1}\}}} \left( \det \widetilde{R}^{(k)}_{i_1\cdots i_{k-1} n, j_1\cdots j_{k-1} n} \right) \left[\sum_{m=1}^{k-1}    (-1)^m    (D_{i_m} h(y'))\times \right.\\
&\left(\sum_{\ell=1}^{m-1}(-1)^{\ell} (D_{j_{\ell}} h(y')) \delta_{i_1j_1}\cdots \delta_{i_{\ell-1}j_{\ell-1}} \delta_{i_{\ell}j_{\ell+1}}\cdots \delta_{i_{m-1}j_{m}}\delta_{i_{m+1}j_{m+1}}\cdots \delta_{i_{k-1}j_{k-1}}\right.\\
&+ (-1)^m (D_{j_m} h(y')) \delta_{i_1j_1}\cdots \delta_{i_{m-1}j_{m-1}} \delta_{i_{m+1}j_{m+1}}\cdots \delta_{i_{k-1}j_{k-1}}+\\
& + \left.\left. \sum_{\ell=m+1}^{k-1} (-1)^{\ell} (D_{j_{\ell}} h(y')) \delta_{i_1j_1}\cdots \delta_{i_{m-1}j_{m-1}} \delta_{i_{m+1}j_{m}}\cdots \delta_{i_{\ell}j_{\ell-1}} \delta_{i_{\ell+1}j_{\ell+1}}\cdots \delta_{i_{k-1}j_{k-1}}      \right)\right].
\end{split}
\end{align}
\end{proposition}
In Section 5 we prove that if $u(x)$ is a strictly $\left(\delta, \widetilde{\gamma}_{k}\right)$-admissible solution to the equation \eqref{1.1}, then $v(y),$ defined by \eqref{1.33}, is a strictly $\left(\widetilde{\delta}, \widetilde{\widetilde{\gamma}}_{k}\right)$-admissible solution to the equation \eqref{1.46} in $\overline{\widetilde{\Omega}_{\rho}}$ with $\widetilde{\delta}=(1+\varepsilon)^{2} \delta,$
$\widetilde{\widetilde{\gamma}}_{k}=\frac{1}{(1+\varepsilon)^{2}} \widetilde{\gamma}_{k},$ where $\varepsilon>0$ is sufficiently small
if $\rho$ is chosen sufficiently small. We denote by $\widetilde{F}_k\left(y', \widetilde{R}\right)$ the corresponding new $k$-Hessian type function of the equation \eqref{1.46}, which is
\begin{equation}\label{1.48}
	\widetilde{F}_{k}(y', \widetilde{R})=\log \left[S_{k}(\widetilde{R})+H_{k}\left(y', \widetilde{R}\right)\right],
\end{equation}
where $H_{k}\left(y^{\prime}, \widetilde{R}\right)$ is defined by \eqref{1.47}. The $\widetilde{d}$-concavity of the function $\widetilde{F}_{k}(y',\widetilde{R})$ will be proved in the following.
\begin{proposition}\label{md1.10}
Suppose $\rho$ is chosen sufficiently small so that
\begin{equation*}
0<\gamma_{k} < \widetilde{\widetilde{\gamma}}_k=\frac{1}{(1+\varepsilon)^2} 	\widetilde{\gamma}_k<1, \quad 0<\widetilde{\delta}= (1+\varepsilon)^2 \delta <\delta_k<1,
\end{equation*}
where  $0<\gamma_{k}< \widetilde{\gamma}_k<1$ are defined in Definition \ref{dn1.3}, $0<\delta<\delta_k<1$ are determined in Theorem \ref{dl1.5}.
 Then for any $y',$ $|y'|\le \rho,$ the function $\widetilde{F}_{k}(y', \widetilde{R})$ is $\widetilde{d}$-concave on the set $D_{\widetilde{\sigma}, \mu(\widetilde{B}),\widetilde{\widetilde{\gamma}}_k}$ in the sense of \eqref{1.21}, where
$$
\mu(\widetilde{B})=\sup_{y \in \overline{\widetilde{\Omega_{\rho}}}\atop z\in \mathbb{R}, p\in \mathbb{R}^n}\|\widetilde{B}(y,z, p)\|
$$
and $C'_2>0$ in \eqref{1.21} does not depend on $y', |y'|<\rho.$
\end{proposition}

 Using \eqref{1.31a} and the  $\widetilde{d}$-concavity of $\widetilde{F}_k(y', \widetilde{R}),$
we show the H\"older continuity of $D^{2} v(y)$ in $ \overline{ \widetilde{\Omega}}_{\rho}$ with some $0<\alpha<1,$ if $v(y)\in C^4(\widetilde{\Omega_{\rho}})\cap C^2(\overline{ \widetilde{\Omega}}_{\rho})$ and $\widetilde{A}(y,z,p), \widetilde{B}(y,z,p), \widetilde{f}(y,z,p)\in C^3(\widetilde{\mathcal{D}}_{\rho}),$
$\widetilde{D}_{\rho}= \overline{ \widetilde{\Omega}}_{\rho} \times \mathbb{R} \times \mathbb{R}^{n}.$ So we will obtain the following theorem at the end of Section 5.

\begin{theorem}\label{dl1.11} Under the assumptions of Theorem \ref{dl1.8} there exist $C_{4}>0,$ $0<\alpha<1,$ that depend on $n, k, \delta, \widetilde{\gamma}_{k}, \beta_{1}, \Omega, A(x,z, p), f(x, z, p), \underline{u}(x), \varphi,$  $\mu_2(B),$ such that if $u(x)$ is any strictly 
$\left(\delta, \widetilde{\gamma}_{k}\right)$-admissible solution to the problem \eqref{1.1}-\eqref{1.2}, the following estimate holds
\begin{equation}\label{1.50}
	\|u\|_{C^{2, \alpha}(\overline{\Omega})} \leq C_{4},
\end{equation}
where $\mu_2(B)$ is defined by \eqref{1.23a}.
\end{theorem}

In Section 6 we study the solvability of the Dirichlet problem \eqref{1.1}-\eqref{1.2} in the classes of strictly admissible solutions.  A necessary condition and some sufficient conditions on $B(x, z, p)$ have been found as follows.
\begin{theorem}[A necessary condition] \label{dl1.12}
	Suppose $0<\delta<1$ and there exists a strictly $\delta$-admissible solution $u(x)$ to the equation \eqref{1.1}, which satisfies the following conditions:
	\begin{itemize}
		\item[(i)\ \ ] $\lambda_u=\inf_{x\in \overline{\Omega}} \lambda_{\min} (\omega(x,u))>0,$
	\item[(ii)\ ] $\mu(B)\le \delta \lambda_u,$
	\item[(iii)] $\sup_{\overline{\Omega}} |u(x)| \le M_0, $  $\sup_{\overline{\Omega}} |Du(x)| \le M_1.$
	\end{itemize}
Then it is necessary that
\begin{equation}\label{1.53}
	\mu(B) \leq\delta\left[\frac{1}{\binom{n}{k} } f_{1}\right]^{\frac{1}{k}},
\end{equation}		
		where $f_1$ is defined by \eqref{1.22a}.
\end{theorem}
The following theorem is the main result of the paper.
\begin{theorem}[Sufficient conditions]\label{dl1.13}
	Suppose $2\le k\le n,$ $0<\delta<1,$ $0<\widetilde{\gamma}_{k} <1$ are defined as in \eqref{1.12a}, $A(x,z,p), f(x,z,p)\in C^3(\mathcal{D}).$ Assume that the following conditions hold:
		\begin{itemize}
		\item[(i)\ \ ] $A(x, z, p) \geq-\gamma_{0}\left(1+|p|^{2}\right) E_{n}, \gamma_{0}>0, $ $ \lambda_{\max }(A(x, z, 0)) \geq 0, $ $D_{z} A(x, z, p) \geq 0;$
		\item[(ii)\ ] $f(x,z,p)>0$ in $\mathcal{D}$ and
		\begin{equation*}
			\inf_{\mathcal{D}}\left[\frac{D_{z} f(x, z,p)}{f(x, z, p)}\right] \geq \frac{k \delta}{(1+\delta^{2})} \beta_{1}, \beta_{1}>0;
				\end{equation*}
			\item[(iii)]  There exists a strictly $\widetilde{\gamma}_{k}$-admissible subsolution
			$\underline{u}(x)\in C^4(\overline{\Omega})$ to the problem
				$$
		S_{k}\left(D^{2} u-A\left(x, u, Du \right)\right) = f\left(x, u, Du\right) \text { in } \Omega, $$
		$$ u=\varphi  \text { on } \partial \Omega,$$
		 that satisfies the following conditions: 
		 $$\lambda_{\underline{u}}>0$$ 
		 and
		 \begin{equation}\label{1.53a}
		 	\gamma_{\underline{u}} > \widetilde{\gamma}_{k}+\varepsilon_{0}, \varepsilon_{0}>0,
		 \end{equation}
		where  $\lambda_{\underline{u}}$ and $\gamma_{\underline{u}}$ are defined by \eqref{1.11}, \eqref{1.22} respectively. Here we assume that $\partial \Omega\in C^4, \varphi\in C^4;$
			\item[(iv)]  Suppose $B(x, z, p) \in B C^{3}(\mathcal{D})$ is a skew-symmetric and satisfies the following conditions:
			\begin{equation}\label{1.55}
				\mu(B)<  \delta \min \left(\lambda_{\underline{u}}, \lambda_{*}\right),
			\end{equation}
		\begin{equation}\label{1.56}
			\mu\left(D_{z} B\right)<\beta_{1} \min \left(\lambda_{\underline{u}}, \lambda_{*}\right),
		\end{equation}
	where
	\begin{equation}\label{1.38a} 
		\lambda_{*}= \widetilde{\gamma}_k\left[\frac{(1+\delta^2)^{-\left[\frac{k}{2}\right]}f_0}{\binom{n}{k} }  \right]^{\frac{1}{k}},
	\end{equation}
	$f_0$ is defined by \eqref{1.22a} with $M_0, M_1$ as in Theorem \ref{dl1.8}. 
	
	Then there exists unique strictly $(\delta, \widetilde{\gamma}_k)$-admissible solution $u(x)$ to the problem \eqref{1.1}-\eqref{1.2} that belongs to $C^{2, \alpha}\left(\overline{\Omega}\right)$ with some $0< \alpha<1,$ where $\alpha$ depends on $n, k, \delta, \widetilde{\gamma}_k,$ $\beta_1,$ $\Omega,$ $A(x,z,p),$ $f(x,z,p),$ $\underline{u}(x),$ $\varphi,$ $\mu_2(B).$
\end{itemize}
\end{theorem}

In the last Section 7, we consider an example of the Dirichlet problem for a
nonsymmetric $k$-Hessian type equation in the cases $2 \leq  k \leq n$ and in the separated case $k=2.$

\section{The comparison principle for the strictly $\delta$-admissible solutions}
First, we prove the following lemma on ellipticity of the equation $\log F_k(R(x,u))=\log f(x,u,Du)$ at a strictly $\delta$-admissible solution.
\begin{lemma}
Suppose $0<\delta<1,$ $\mu>0$ and  $R=\omega+\beta \in D_{\delta, \mu}.$
Then for $F_{k}(R)=\log \left(S_{k}(R)\right)$ we have
\begin{equation}\label{2.1}
\frac{k}{n} (1+\delta^2)^{-2\left[\frac{k}{2}\right]} \frac{\lambda_{\min}^k(\omega)}{\lambda_{\max}^{k+1}(\omega) } |\xi|^2
\le \frac{1}{2} \sum_{i, j=1}^{n}\left(\frac{\partial F_k(R)}{\partial R_{ij}} + \frac{\partial F_k(R)}{\partial R_{ji}}  \right)\xi_i \xi_j  	\le \frac{ (1+\delta^2)^{\left[\frac{k}{2}\right]}}{\lambda_{\min}(\omega)}|\xi|^2	\end{equation}
for any $\xi=(\xi_1,\cdots,\xi_n)^T\in \mathbb{R}^n.$
\end{lemma}
\begin{proof}
	Suppose
	$$R=\omega+\beta=C^{-1}(D+C\beta C^{-1})C= C^{-1}(D+\widetilde{\beta})C=C^{-1} \widetilde{R} C=\left[R_{ij}\right]_{n\times n},$$ 
	where $C$ is an orthogonal matrix, $D= \operatorname{diag}\left(\lambda_{1}, \ldots, \lambda_{n}\right), \quad \lambda_{j}>0.$
	Since $S_{k}(R)=S_{k}(\widetilde{R}),$
	then $F_k(R)=F_{k}(\widetilde{R}).$ We denote
	$$
	\eta=\left(\eta_{1}, \cdots, \eta_{n}\right)^{T}=C \xi.
	$$
	Then we have
	$$\frac{1}{2} \sum_{i, j=1}^{n}\left(\frac{\partial F_k(R)}{\partial R_{ij}} + \frac{\partial F_k(R)}{\partial R_{ji}}  \right)\xi_i \xi_j=
	\frac{1}{2} \sum_{i, j=1}^{n}\left(\frac{\partial F_k(\widetilde{R})}{\partial \widetilde{R}_{ij}} + \frac{\partial F_k(\widetilde{R})}{\partial \widetilde{R}_{ji}}  \right)\eta_i \eta_j. $$
	So, we can assume that $R=D+\beta \in D_{\delta, \mu}.$
	We note that if $\sigma=D^{-\frac{1}{2}} \beta D^{-\frac{1}{2}},$ then $\|\sigma\| \leq \delta.$
	
	To prove \eqref{2.1} we recall now some facts
	from \cite{2}. If for indices $i_{1} i_{2} \ldots i_{k}$ with
	$1 \leq  i_{1}<\cdots<i_{k}   \leq n$ we set
	$$
	R_{i_{1} \ldots i_{k}}=\left[R_{i_pi_q}\right]_{p, q=1}^{k}, G_{i_{1} \ldots i_{k}}(R)=\det\left(R_{i_{1} \ldots i_{k}}\right),
	\left(R_{i_{1}} \ldots i_{k}\right)^{-1}=\left[\left(R_{i_{1}} \ldots i_{k}\right)_{i_{p} i_q}^{-1}\right]_{p, q=1}^{k},$$
	then we have 
	\begin{equation}\label{2.2}
		\frac{\partial F_k(R)}{\partial R_{ij}}= \frac{1}{S_k(R)} \sum_{1\le i_1<\cdots <i_k\le n} G_{i_1\cdots i_k}(R) \sum_{p,q=1 }^k(R_{i_{1} \ldots i_{k}})^{-1}_{i_q i_p} \delta_{i i_p}\delta_{j i_q}.
	\end{equation}
There are some following relations:
\begin{equation*}
	\frac{\left(R_{i_{1}\cdots i_{k}}\right)^{-1}+\left[\left(R_{i_{1} \ldots i_{k}}\right)^{-1}\right]^{T}}{2}=D_{i_1\cdots i_k}^{-\frac{1}{2}}\left(E_{i_1\cdots i_k} -\sigma_{i_1\cdots i_k}^{2} \right) D_{i_1\cdots i_k}^{-\frac{1}{2}},
\end{equation*}
\begin{equation}\label{2.2b}
(1+\delta^2)^{-\left[\frac{k}{2}\right]} E_{i_{1}\cdots i_{k}}\le \left( E_{i_{1}\cdots i_{k}}- \sigma_{i_1\cdots i_k}^{2}\right)^{-1}\le  E_{i_{1}\cdots i_{k}},	
\end{equation}
where 
$$E_{i_{1}\cdots i_{k}}= \left[\delta_{i_p i_q}\right]^k_{p,q=1}, D_{i_1\cdots i_k}^{-\frac{1}{2}}=\operatorname{diag}\left( \lambda_{i_{1}}^{-\frac{1}{2}},\cdots, \lambda_{i_{k}}^{-\frac{1}{2}} \right),$$
$$\sigma_{i_{1} \cdots i_{k}}=D_{i_{1}\cdots i_{k}}^{-\frac{1}{2}} \beta_{i_{1} \ldots i_{k}} D_{i_{1}\cdots i_{k}}^{-\frac{1}{2}},\left\|\sigma_{i_{1}\cdots i_{k}}\right\| \leq  \delta,$$
\begin{equation}\label{2.2c}
	G_{i_1\cdots i_k} (D)\le G_{i_1\cdots i_k} (R)\le (1+\delta^2)^{\left[\frac{k}{2}\right]} G_{i_1\cdots i_k} (D),
\end{equation}
\begin{equation}\label{2.5}
	S_k(D)\le S_k(R)\le (1+\delta^2)^{\left[\frac{k}{2}\right]} S_k(D),
\end{equation}
\begin{equation}\label{2.6}
	(1+\delta^2)^{-\left[\frac{k}{2}\right]} \frac{G_{i_1\cdots i_k} (D) }{S_k(D)} \le \frac{G_{i_1\cdots i_k} (R) }{S_k(R)} \le (1+\delta^2)^{\left[\frac{k}{2}\right]} \frac{G_{i_1\cdots i_k} (D) }{S_k(D)}.
\end{equation}
From \eqref{2.2}-\eqref{2.6} it follows that
\begin{equation}\label{2.7}
	\begin{split}	\frac{1}{2} \sum_{i, j=1}^{n}\left(\frac{\partial F_{k}(R)}{\partial R_{ij}}+\frac{\partial F_{k}(R)}{\partial R_{j i}}\right) \xi_{i} \xi_{j} &\leq 
	\frac{\left(1+\delta^{2}\right)^{\left[\frac{k}{2}\right]}}{\lambda_{\min }} \sum_{1 \leq  i_{1} <\cdots<i_{k}\leq  n}   \frac{\lambda_{i_{1}} \cdots\lambda_{i_k}}{\sigma_{k}(\lambda)}\left(\sum_{p=1}^{k} \xi_{i_{p}}^{2}\right)\\
	&\le \frac{\left(1+\delta^{2}\right)^{\left[\frac{k}{2}\right]}}{\lambda_{\min} }|\xi|^{2}.
	\end{split}
	\end{equation}
Here we have used the facts that $\sum_{p=1}^{k} \xi_{i_{p}}^{2} \leq |\xi|^{2}$ and $\frac{1}{\sigma_{k}(\lambda)} \sum_{1 \leq i_{1}<\cdots i_{k} \leq n} \lambda_{i_1} \cdots \lambda_{i_k}=1.$

On other side, we also have from \eqref{2.2}-\eqref{2.6} that 

\begin{equation}\label{2.8}
\begin{split}
	\frac{1}{2} \sum_{i, j=1}^{n}\Big(\frac{\partial F_{k}(R)}{\partial R_{ij}}+&\frac{\partial F_{k}(R)}{\partial R_{j i}}\Big) \xi_{i} \xi_{j}\\
	&\ge \frac{\left(1+\delta^{2}\right)^{-2\left[\frac{k}{2}\right]}}{\lambda_{\max} }   
	\sum_{1 \leq  i_{1} <\cdots<i_{k}\leq  n}   \frac{\lambda_{i_{1}} \cdots\lambda_{i_k}}{\sigma_{k}(\lambda)} \sum_{p=1}^{n} \xi_{i_{p}}^{2} \\
	&= \frac{\left(1+\delta^{2}\right)^{-2\left[\frac{k}{2}\right]}}{\lambda_{\max} } \sum_{\ell=1}^{n} \left[ \sum_{1 \leq  i_{1} <\cdots<i_{k}\leq  n\atop \ell\in \{i_1,\cdots,i_k\}}   \frac{\lambda_{i_{1}} \cdots\lambda_{i_k}}{\sigma_{k}(\lambda)}  \right]\xi^2_{\ell}\\
	&=\frac{\left(1+\delta^{2}\right)^{-2\left[\frac{k}{2}\right]}}{\lambda_{\max} } \sum_{\ell=1}^{n} \frac{ \lambda_{\ell} \sigma_{k-1}^{(\ell)}(\lambda)}{\sigma_{k}(\lambda)}\xi^2_{\ell}\\
	&\ge \frac{\left(1+\delta^{2}\right)^{-2\left[\frac{k}{2}\right]}}{\lambda_{\max}} \frac{\binom{n-1}{k-1}\lambda_{\min}^k}{\binom{n}{k} \lambda_{\max}^k} |\xi|^2\\
	&=\frac{k}{n} \left(1+\delta^{2}\right)^{-2\left[\frac{k}{2}\right]} \frac{ \lambda_{\min}^k}{\lambda_{\max}^{k+1}}|\xi|^2.
\end{split}
\end{equation}
Then \eqref{2.1} follows from \eqref{2.7} and \eqref{2.8}.
\end{proof}
For $u(x) \in C^{2}(\overline{\Omega})$ we set
$$
G_{k}[u](x)=\log \left(S_{k}(R(x, u))\right)-\log f(x, u, D u),
$$
where $R\left(x, u\right)$ is defined by \eqref{1.7a}.

\begin{theorem}\label{dl2.2}
	Suppose $A(x,z,p), f(x,z,p)\in C^1(\mathcal{D}),$ $B(x,z,p)\in BC^1(\mathcal{D}),$ $0<\delta<1$ and	suppose $u(x), v(x) \in C^{2}(\overline{\Omega})$
	and satisfy the following conditions
	\begin{itemize}
		\item[(i)\ \ ] $G_{k}[u](x) \leq G_{k}[v](x),\ x \in \Omega;$
			\item[(ii)\ ] $\lambda_u>0,\ \lambda_v>0;$
				\item[(iii)] $D_z A(x,z,p)\ge 0,\ (x,z,p)\in \mathcal{D};$
				\item[(iv)\ ]	$\mu(B) \leq  \delta \min \left(\lambda_{u}, \lambda_{v}\right);$
			\item[(v)\ \ ]			$\mu\left(D_{z}B\right) \leq \beta_{1} \min \left(\lambda_{u}, \lambda_v\right),\ \beta_{1}>0;$
			\item[(vi)\ ]	$f(x, z, p)>0,\ (x, z, p) \in   \mathcal{D};$
				\item[(vii)]		$\inf_{\mathcal{D}}\left[\frac{D_{z} f\left(x, z, p\right)}{f(x, z, p)}\right] \geq  \frac{k \delta}{(1+\delta^{2})} \beta_{1}.$
			\end{itemize}	
Then the following assertions are true:
\begin{itemize}
	\item[(a)] If $u(x) \geq  v(x)$ on $\partial \Omega,$ then
	$$
	u(x) \geq  v(x) \text { in } \Omega,  
	$$
		\item[(b)] If $u(x)=v(x)$ on $\partial \Omega,$ then
		$$
		\frac{\partial u(x)}{\partial \nu} \geqslant \frac{\partial v(x)}{\partial \nu} \text { on } \partial \Omega,
		$$
		where $\nu$ is the unit inward normal
		at $x \in \partial \Omega.$
		\end{itemize}
	\end{theorem}
\begin{proof}
From the assumptions (ii) and (iv) it follows that for any $x \in \Omega$
$$
R(x, u), R(x, v) \in D_{\delta, \mu(B)}
$$
and
$$
\lambda_{\min }(\omega(x, u)) \geq\lambda_{u}>0,\
\lambda_{\min }\left(\omega\left(x, v\right)\right) \geq \lambda_{v}>0.$$
Then, by using \eqref{2.1}-\eqref{2.6}   and the following
relation (\cite{11})
\begin{equation}\label{2.7N}
	\frac{\left(R_{i_{1} \cdots i_{k}}\right)^{-1}-\left[\left(R_{i_1 \cdots  i_{k}}\right)^{-1}\right]^{T}}{2}= D_{i_{1} \cdots i_{k}}^{-\frac{1}{2}} (-\sigma_{i_{1} \cdots i_{k}})(E_{i_{1} \cdots i_{k}}-\sigma^2_{i_{1} \cdots i_{k}}) D_{i_{1} \cdots i_{k}}^{-\frac{1}{2}},
\end{equation}
the proof of the theorem will go analogously as in the proof of the comparison principle for nonsymmetric Monge-Ampère type equations \eqref{1.10} in \cite{11} (Theorem 4).
\end{proof}
\section{The $C^{2}(\overline{\Omega})$-estimates for strictly $\left(\delta, \widetilde{\gamma}_{k}\right)$-admissible solutions and the H\"older continuity of their second-order derivatives inside the domain}
We recall that for $u(x)\in C^{2}(\overline{\Omega}), x\in \overline{\Omega}$ the matrices $\omega(x,u)$ and $R(x,u)$ are defined respectively by \eqref{1.7} and \eqref{1.7a}. 

The equation \eqref{1.1} can be written as 
\begin{equation}\label{3.1}
S_k(R(x,u))=f(x,u,Du),\ x \in \overline{\Omega}.
\end{equation}
\subsection{Proof of Proposition \ref{md1.7}}
\begin{itemize}
\item[(i)] Suppose $u(x)$ is a strictly admissible solution to \eqref{3.1}, i.e. $\lambda_{\min }(\omega(x,u))\ge \lambda_u>0.$
We have from \eqref{2.5} and \eqref{3.1} that
\begin{align*}
	\binom{n}{k} \lambda_{\min }^k(\omega(x,u))&\le S_k(\omega(x,u))\le S_k(R(x,u))=f(x,u,Du)\\
	&\le 	\sup_{x\in \overline{\Omega}\atop |z|\le M_0, |p|\le M_1} \left[f(x,z,p)\right]=f_1,
\end{align*}
if $|u(x)|\le M_0, |Du(x)|\le M_1.$ So, \eqref{1.23} is proved.
\item[(ii)]  Suppose $u(x)$ is a strictly $\delta$-admissible solution to \eqref{3.1}.  Then 
$$
\mu(B) \leq  \delta \lambda_{u}
$$
and $R(x, u) \in D_{\delta, \mu(B)}$ for any $x \in \overline{\Omega}.$
From \eqref{3.1} and \eqref{2.5}  we obtain
\begin{align*}
	\binom{n}{k}\left(1+\delta^{2}\right)^{\left[\frac{k}{2}\right]} \lambda_{\max}^k(\omega(x,u))&\ge \left(1+\delta^{2}\right)^{\left[\frac{k}{2}\right]} S_k (\omega(x,u))\\
	&\ge S_k (R(x,u))=f(x,u,Du)\ge f_0,
\end{align*}
from which it follows \eqref{1.24}.
\item[(iii)] If $u(x)$ is a strictly $\widetilde{\gamma}_k$-admissible solution to \eqref{3.1}, then from \eqref{1.22} implies that
$$\lambda_{\min }(\omega(x,u))\ge \widetilde{\gamma}_k \lambda_{\max }(\omega(x,u)),\ x\in \overline{\Omega} $$
and \eqref{1.25a} follows therefore from \eqref{1.23} and the last inequality.
\item[(iv)] If $u(x)$ is a strictly $(\delta,\widetilde{\gamma}_k)$-admissible solution to \eqref{3.1}, then  \eqref{1.25} and  \eqref{1.26} follow from \eqref{1.23}, \eqref{1.24}, \eqref{1.25a} and the last inequality.
\end{itemize}
\hfill $\square$
\subsection{Proof of Theorem \ref{dl1.8}}
Suppose $u(x)$ is a strictly $\left(\delta, \widetilde{\gamma}_k\right)$-admissible solution to the problem \eqref{1.1}-\eqref{1.2}, it is also a strictly $\delta$-admissible  one. Since there exists a strictly $\widetilde{\gamma}_k$-admissible subsolution $\underline{u}(x)$ to \eqref{1.27}-\eqref{1.28}, this function due to \eqref{2.5} is also a strictly  $\widetilde{\gamma}_k$-admissible subsolution to the problem \eqref{1.1}-\eqref{1.2}. From the condition $(v)$ it follows that the function $\underline{u}(x)$ is also a strictly $\delta$-admissible solution to \eqref{1.1}-\eqref{1.2}. Therefore we can apply the comparison principle (Theorem \ref{dl2.2}) for $u(x)$ and $\underline{u}(x)$ to conclude that 
$u \geq \underline{u}$ in $\Omega, \frac{\partial u}{\partial \nu} \geq \frac{\partial \underline{u}}{\partial \nu}$ on $\partial \Omega,$
where $\nu$ is the unit inner normal to $\partial \Omega.$

By using this fact and by following the same arguments as in \cite{8}, from the structure conditions for $A(x, z, p),$ we can obtain the following estimates
\begin{equation*}
\sup_{\overline{\Omega}}|u| \leq  M_{0}, \quad 	\sup_{\overline{\Omega}}|Du| \leq  M_{1},
\end{equation*}
where $M_{0}$ depends on $|\underline{u}|_{0, \overline{\Omega}}, |\varphi|_{0,\overline{\Omega}}$ and $M_{1}$ depends on $n, \gamma_{0},|\underline{u}|_{1,\overline{\Omega}},  |\varphi|_{2, \overline{\Omega}}$ and $\Omega.$
We prove that there exists $M_{2}>0$ such that
\begin{equation}\label{3.3}
\sup_{\Omega}|D^2u| \leq  M_{2}.
\end{equation}

Indeed, since $u(x)$ is a strictly $\left(\delta, \widetilde{\gamma}_{k}\right)$-admissible solution, then it follows from \eqref{1.26} that for any $x \in \overline{\Omega}$ we have
\begin{equation}\label{3.4}
\left|\omega\left(x, u\right)\right| \leq \sqrt{n}\|\omega(x, u)\|=\sqrt{n} \lambda_{\max}(\omega(x, u))
\le\frac{\sqrt{n}}{\widetilde{\gamma}_k}\left[\frac{1}{ \binom{n}{k}} f_{1}\right]^{\frac{1}{k}}.
\end{equation}
From the equality
$$
D^{2} u=\omega(x, u)+A\left(x, u(x), D u(x)\right), \ x\in \overline{\Omega}
$$
and from \eqref{3.4}, we obtain \eqref{3.3}, where $M_{2}$ depends on $n, k, \widetilde{\gamma}_{k}, M_{0}, M_{1},$ $A(x, z, p)$ and $f(x, z, p).$

\hfill $\square$
\subsection{H\"older continuity of the second-oder derivatives inside the domain}
\begin{proposition}\label{md3.3}
Suppose $A(x, z, p),$ $f(x, z, p)\in C^{3}(\mathcal{D}),$ $ B(x, z, p) \in B C^{3}(\mathcal{D}), u(x) \in C^{4}(\Omega)$ is a strictly $\left(\delta, \widetilde{\gamma}_k \right)$-admissible solution to the problem \eqref{1.1}-\eqref{1.2}. Then for any $\Omega' \subset \subset  \Omega$
there exist $C_{4}^{'}>0,0<\alpha<1$ such that
\begin{equation}\label{3.3a}
	\left\|D^{2} u\right\|_{C^{2, \alpha} \left(\overline{\Omega}'\right)}\leq C_{4}^{'},
\end{equation}
where $C'_4$ and $\alpha$ depend on $n, k, \delta, \widetilde{\gamma}_k,$ $\beta_1,$ $\Omega',$ $A,$ $f,$ $\underline{u},$  $\mu_2(B).$
\end{proposition}
First we prove the following lemma on uniform ellipticity of the equation $\log\left(F_k(R(x,u))\right)=\log f(x,u,Du)$ at a $(\delta, \widetilde{\gamma}_k)$-admissble solution $u(x)$ by improving \eqref{2.1}.
\begin{lemma}
Suppose $u(x)$ is a strictly $\left(\delta, \widetilde{\gamma}_{k}\right)$-admissible solution to the problem \eqref{1.1}-\eqref{1.2}.
Then for $R=R(x, u)=\left[R_{ij}(x, u)\right]$ the
following estimates are true for $x\in \overline{\Omega
}:$
\begin{equation}\label{3.6}
	\frac{k {\binom{n}{k}}^{\frac{2}{k}} \widetilde{\gamma}_{k}^{k-2}  \left(1+\delta^{2}\right)^{-2\left[\frac{k}{2}\right]}}{n f_{1}^{\frac{1}{k}}}|\xi|^2 \le 
	\frac{1}{2} \sum_{i, j=1}^{n}\left(\frac{\partial F_{k}(R)}{\partial R_{ij}}+\frac{\partial F_{k}(R)}{\partial R_{ji}}\right) \xi_i\xi_j\le
	\frac{{\binom{n}{k}}^{\frac{2}{k}}\left(1+\delta^{2}\right)^{3\left[\frac{k}{2}\right]}}{\widetilde{\gamma}_k^{2} f_{0}^{\frac{1}{k}}}|\xi|^{2},
\end{equation}	
where $f_0, f_1$ are defined by \eqref{1.22a}.
\end{lemma}
\begin{proof}
Since $R\in D_{\delta,\mu(B)},$ the inequalities \eqref{2.1} are true. Then \eqref{3.6} follows from \eqref{2.1}, the relation 
$$
\lambda_{\min }(\omega(x, u)) \geq \widetilde{\gamma}_k  \lambda_{\max }\left(\omega\left(x, u\right)\right),\ x\in \overline{\Omega}
$$
and from \eqref{1.25}, \eqref{1.26}.
\end{proof}
\begin{proof}[Proof of Proposition \ref{md3.3}]
To prove the H\"older continuity inside $\Omega$ for second-order derivatives of the solution $u(x)$ we consider the equation $F_k(R(x,u))=\log f(x,u,Du)$ in $\Omega$ and we can use the following already established facts:
\begin{itemize}
	\item[(i)\ \ ] The $C^{2}(\overline{\Omega})$-estimates \eqref{1.29} 
	$$
	\|u\|_{C^{2}(\overline{\Omega})} \leq  C_{3};
	$$
	\item[(ii)\ ] The uniform ellipticity \eqref{3.6} of the equation \eqref{1.1} at $R=R(x, u)$ for any $x \in \overline{\Omega};$
	\item[(iii)] The strict concavity \eqref{1.19} of 
	$F_{k}(\omega+\beta)=\log \left(S_{k}(\omega+\beta)\right)$ as a function
	of $\omega>0$ when $\beta^T=-\beta$ is fixed, i.e 
	$$
	d^{2} F_{k}(R, P) \leq -\frac{C_{1}}{\lambda_{\max}^{2}(\omega)}.|P|^{2},\ P^{T}=P,\ C_{1}>0,
	$$
	where $\lambda_{\max}(\omega(x,u))$ satisfies the estimates \eqref{1.26};
	\item[(iv)] The following version \eqref{1.21} of the $d$-concavity of the function $F_k(R)$ on the set $D_{\delta,\mu(B),\widetilde{\gamma}_k}:$
	$$F_{k}\left(R^{(1)}\right)-F_{k}\left(R^{(0)}\right) \leq \sum_{i, j=1}^{n} \frac{\partial F_{k}\left(R^{(0)}\right)}{\partial R_{ij}}\left(R_{i j}^{(1)}-R_{i j}^{(0)}\right)
	+C_{2} \frac{\left|\beta^{(1)}-\beta^{(0)}\right|^{2}}{\lambda_{\min }^{2}\left(\omega^{(\tau)}\right)},$$
	where $R^{(0)}=\omega^{(0)}+\beta^{(0)}, R^{(1)}=\omega^{(1)}+\beta^{(1)}\in D_{\delta,\mu(B),\widetilde{\gamma}_k},$ 
	$\omega^{(\tau)}=(1-\tau)\omega^{(0)}+\tau\omega^{(1)}, 0<\tau<1, \lambda_{\min}(\omega(x))$ satisfies the estimates \eqref{1.25}. 
	
	Hence, the facts mentioned above and the methods of L.C. Evans and N.V. Krylov allow ones with the aid of \eqref{2.2}-\eqref{2.7N} to get the desired H\"older continuity \eqref{3.3a} of $D^2u$ inside $\Omega$ (see \cite{10}, Section 17.4).
\end{itemize} 
\end{proof}
\section{A new kind of the $k$-Hessian type equation in a neighborhood of the boundary}
\subsection{The $k$-compound of a square matrix}
Let $M=\left[M_{ij}\right]$ be an $n \times n$ matrix with entries in $\mathbb{R}$ or $\mathbb{C}.$ Suppose that $i_{1} i_{2}  \cdots i_k$
and $j_{1} j_{2}\cdots j_{k}$
are indices such that
$$
1 \leqslant i_{1}<\cdots<i_{k} \leqslant n, \quad 1 \leqslant j_{1}<\cdots <j_{k} \leqslant n.
$$
We denote
\begin{equation*}
M^{(k)}_{i_{1}\cdots i_{k},j_{1}\cdots j_{k}}=\left[M_{i_pj_q}\right]^{k}_{p,q=1}.
\end{equation*}
Then $\det\left(	M^{(k)}_{i_{1}\cdots i_{k},j_{1}\cdots j_{k}} \right) $
is a minor at the intersection of the rows $i_{1}, i_{2}, \cdots, i_{k}$ and the columns $j_{1}, j_{2}, \cdots, j_{k}.$ When the 
indices $i_{1} i_{2} \cdots i_{k}$ are arranged in the lexical order, the resulting $\binom{n}{k}\times \binom{n}{k} $ square matrix, that consists of corresponding minors,  is called the $k$-compound of the matrix $M$
and written as $M^{(k)}.$ That means
\begin{equation*}
M^{(k)}=\left[ \det \left(  	M^{(k)}_{i_{1}\cdots i_{k},j_{1}\cdots j_{k}}   \right)   \right]_{ \binom{n}{k} \times \binom{n}{k}}.
\end{equation*}
We list here some properties of the $k$-compounds.
\begin{proposition}[\cite{1}]\label{md4.1}
Let $M$ and $N$ be matrices in $\mathbb{C}^{n \times n}.$ Then the following assertions are true:
\begin{itemize}
	\item[(i)\ \ \ ] Binet-Cauchy Theorem
	$$
	(M N)^{(k)}=M^{(k)} N^{(k)};
	$$
	\item[(ii)\ \ ] $\left(M^{(k)}\right)^{T}=\left(M^{T}\right)^{(k)};$
	\item[(iii)\ ]
	$\overline{M^{(k)}}=(\overline{M})^{(k)};$
	\item[(iv)\ \ ]
	$(M^{(k)})^{*}=\left(M^{*}\right)^{(k)}, \quad M^{*}=(\overline{M})^{T};$	
	\item[(v)\ \ \ ]	
	$M$ is non-singular if and only if $M^{(k)}$ is non-singular, and
	$$
	[M^{(k)}]^{-1}=\left(M^{-1}\right)^{(k)};
	$$
	\item[(vi)\ \ ]	
	$(h M)^{(k)}=h^{k} M^{(k)},$ for any $h \in \mathbb{C};$
	\item[(vii)\ ]		$M^{(k)}$ is symmetric if $M$ is symmetric;
	\item[(viii)]	
	If $M=\operatorname{diag}\left(\lambda_{1}, \lambda_{2}, \cdots, \lambda_{n}\right) \in \mathbb{C}^{n \times n},$ then
	$$M^{(k)}=\operatorname{diag}\left(\lambda_{i_{1}} \lambda_{i_{2}} \cdots \lambda_{i_{k}}; 1 \leqslant i_{1}<\ldots<i_{k}  \leqslant n\right).$$
	\item[(ix)\ \ ] If $M\in \mathbb{C}^{n \times n},$ then 
	$$ S_k(M)=\sigma_k(\lambda(M))=\operatorname{Tr} \left(M^{(k)}\right)=
	\sum_{1 \leq  i_{1} <\cdots<i_{k}\leq  n} \det   \left(  	M^{(k)}_{i_{1}\cdots i_{k},i_{1}\cdots i_{k}}   \right) .$$
\end{itemize}
\end{proposition}

\subsection{Proof of Proposition \ref{md1.9}}
By using Proposition \ref{md4.1} we rewrite the left hand side of \eqref{1.41}  as follows
\begin{equation}\label{4.3}
\begin{split}
	S_{k} \Big[ J^T(D^{2} v -&\widetilde{A}(y, v, D v)-\widetilde{B}(y, v, D v)\big) J \Big]\\
	&=S_k\left( J^T \widetilde{R} J \right)=\operatorname{Tr} \left( {\left( J^T \widetilde{R} J \right)}^{(k)}   \right) \\
	&=\operatorname{Tr} \left( (J^T)^{(k)}  (\widetilde{R})^{(k)} J^{(k)}  \right)
	=\operatorname{Tr} \left(J^{(k)} (J^T)^{(k)}   (\widetilde{R})^{(k)}    \right)\\
	& =\operatorname{Tr} \left( (J J^T)^{(k)} \widetilde{R}^{(k)}     \right)\\
	& =\sum_{1\le i_1<\cdots <i_{k}\le n\atop 1\le j_1<\cdots <j_{k}\le n } \det \left( (J J^T)^{(k)}_{i_i\cdots i_k, j_1\cdots j_k}\right)\det \left( \widetilde{R}^{(k)}_{i_i\cdots i_k, j_1\cdots j_k}\right) 
\end{split}
\end{equation}
From \eqref{1.35} it follows 
\begin{equation}\label{4.4}
JJ^T=\left[\begin{array}{cc}
	E_{n-1} &  -(D h)^{T}\\
	 -D h& 1+|D h|^{2}
\end{array}\right].
\end{equation}
Then \eqref{1.46} and the proposition \ref{md1.9} follow from \eqref{4.3} and the following lemma.
\hfill $\square$
\begin{lemma}
The entries of $(JJ^T)^{(k)}$ are of the following values:
\begin{itemize}
	\item[(i)\ \ ] If $1\le i_1< \cdots < i_k\le n-1,$ $ 1\le j_1< \cdots < j_k\le n-1$  then 
	$$\det (JJ^T)^{(k)}_{i_i\cdots i_k, j_1\cdots j_k}
	=\delta_{i_1j_1}  \delta_{i_2j_2}\cdots \delta_{i_kj_k}; $$ 
	\item[(ii)\ ] If $1\le i_1< \cdots < i_{k-1}\le n-1,$   then 
	$$\det (JJ^T)^{(k)}_{i_1\cdots i_{k-1}n, i_1\cdots i_{k-1}n}
	=\left( 1+\sum_{m\in \{1,\cdots,n-1\}\backslash\{i_1,\cdots,i_{k-1}\}} (D_m h)^2    \right); $$
	\item[(iii)] If $1\le i_1< \cdots < i_k\le n-1,$ $ 1\le j_1< \cdots < j_{k-1}< n=j_k$  then 
	\begin{align*}
		\det (JJ^T)^{(k)}_{i_i\cdots i_k, j_1\cdots j_{k-1} n}
		&=\det (JJ^T)^{(k)}_{j_i\cdots j_{k-1} n, i_1\cdots i_{k} }\\
		&=(-1)^{k-1} \sum_{m=1}^{k} (-1)^m (D_{i_m}h) \delta_{i_1j_1} \cdots  \delta_{i_{m-1}j_{m-1}}\delta_{i_{m+1}j_{m}}   \cdots \delta_{i_kj_{k-1}}	;
	\end{align*}

	\item[(iv)\ ] If $1\le i_1< \cdots < i_{k-1}\le n-1,$ $ 1\le j_1< \cdots < j_{k-1}\le n-1$ with $(i_1,\cdots, i_{k-1})\ne (j_1,\cdots, j_{k-1}),$ then 
	\begin{align*}
		& \det (JJ^T)^{(k)}_{i_1\cdots i_{k-1}n, j_1\cdots j_{k-1}n}
		=\det (JJ^T)^{(k)}_{j_1\cdots j_{k-1}n, i_1\cdots i_{k-1}n}\\
		&= - \sum_{m=1}^{k-1} (D_{i_m}h)  (-1)^m  \left[ \sum_{\ell=1}^{m-1}(-1)^{\ell} (D_{j_{\ell}}h)  \delta_{i_1j_1} \cdots  \delta_{i_{\ell-1}j_{\ell-1}}\delta_{i_{\ell}j_{\ell+1}}   \cdots \delta_{i_{m-1}j_{m}}  \delta_{i_{m+1}j_{m+1}} \cdots \delta_{i_{k-1}j_{k-1}}   \right.\\
		&+ (-1)^m (D_{j_m}h)     \delta_{i_1j_1} \cdots  \delta_{i_{m-1}j_{m-1}} \delta_{i_{m+1}j_{m+1}}   \cdots \delta_{i_{k-1}j_{k-1}}  \\
		&+\left. \sum_{\ell=m+1}^{k-1}(-1)^{\ell} (D_{j_{\ell}}h) 
		\delta_{i_1j_1} \cdots  \delta_{i_{m-1}j_{m-1}}\delta_{i_{m+1}j_{m}}   \cdots \delta_{i_{\ell}j_{\ell-1}}  \delta_{i_{\ell+1}j_{\ell+1}} \cdots \delta_{i_{k-1}j_{k-1}}
		\right].
	\end{align*}
\end{itemize}
\end{lemma}
\begin{proof}
First, we prove that if $1\le \ell\le n,$ $1\le i'_1<i'_2<\cdots < i'_{\ell}\le n,$
$1\le j'_1<j'_2<\cdots < j'_{\ell}\le n,$ then 
\begin{equation}\label{4.5}
	\det \begin{bmatrix}
		\delta_{i'_{1}j'_{1}} & \delta_{i'_{1}j'_{2}}& \cdots 	&\delta_{i'_{1}j'_{\ell}}					\\
		\delta_{i'_{2}j'_{1}} & \delta_{i'_{2}j'_{2}}& \cdots 	&\delta_{i'_{2}j'_{\ell}}					\\
		\cdots & 		\cdots &  \cdots &\cdots \\
		\delta_{i'_{\ell}j'_{1}} & \delta_{i'_{\ell}j'_{2}}& \cdots 	&\delta_{i'_{\ell}j'_{\ell}}					\\
	\end{bmatrix}
	= \delta_{i'_{1}j'_{1}} \cdots \delta_{i'_{\ell}j'_{\ell}}	.
\end{equation}
Indeed, the determinant is not zero if and only if all the following conditions hold: there exists $j'_{m_1},$ $1\le m_1\le \ell ,$ such that $j'_{m_1}=i'_1,$ there exists $j'_{m_2},$ $1\le m_2\le \ell ,$ $m_2> m_1$ such that $ j'_{m_2}=i'_2, \cdots,$ there exists $j'_{m_{\ell}},$ $1\le m_{\ell}\le \ell ,$ $m_{\ell}> m_p,$ $1\le p\le \ell-1$ such that $j'_{m_{\ell}}=i'_{\ell}.$ But these conditions hold if and only if 
$$\{  i'_1, i'_2, \cdots, i'_{\ell} \}=\{  j'_1, j'_2, \cdots, j'_{\ell} \}.$$
By using \eqref{4.5}, the entries of the matrix $(JJ^T)^{(k)}$ can be calculated directly from \eqref{4.4}.
\end{proof}
\section{The $\widetilde{d}$-concavity of the new kind of the $k$-Hessian type function and the $C^{2,\alpha}$ estimates}
\subsection{Proof of Proposition \ref{md1.10}}
Suppose $J$ and $J^T$ are defined by \eqref{1.35}. We set 
\begin{equation*}
S=JJ^T=
\left[\begin{array}{cc}E_{n-1} & -D h \\ -(D h)^{T} & 1+|D h|^{2}\end{array}\right].
\end{equation*}
We denote the eigenvalues of $S$ as $s_{1}, \cdots, s_{n}$
with $s_{1} \geqslant s_{2} \geqslant \cdots \geqslant s_{n}.$ One can verify
that $s_{2}=s_{3}=\cdots=s_{n-1}=1$ and
\begin{align*}
s_{1}&=\frac{\left(2+|D h|^{2}\right)+\sqrt{\left(2+|D h|^{2}\right)^{2}-4}}{2} \geqslant 1,\\
s_{n}&=\frac{\left(2+|D h|^{2}\right)-\sqrt{\left(2+\left|D h\right|^{2}\right)^{2}-4}}{2}=\frac{1}{s_{1}}.
\end{align*}
We have 
$$\frac{1}{\sqrt{s_{1}}} \leqslant\|J\|=\|J^T\|=\sqrt{\left\|J^{T} J\right\|} = \sqrt{s_{1}}.$$
But $s_1, \cdots, s_n$ are also the eigenvalues of the matrix  $S^{-1}=\left(J^{T}\right)^{-1}J^{-1}.$ So we have 
$$\frac{1}{\sqrt{s_{1}}} \leqslant\|J^{-1}\|=\|(J^{T} )^{-1}\|= \sqrt{s_{1}}.$$
Therefore, we can assume that the neighborhood $\widetilde{\Omega}_{\rho}$ is chosen sufficiently small so that
\begin{equation}\label{5.2}
\frac{1}{\sqrt{1+\varepsilon}} \leqslant\|J\|=\left\|J^{T}\right\| \leqslant \sqrt{1+\varepsilon}, 
\end{equation}
\begin{equation}\label{5.3}
\frac{1}{\sqrt{1+\varepsilon}} \leqslant\left\|J^{-1}\right\|=\left\|\left(J^{T}\right)^{-1}\right\| \leqslant \sqrt{1+\varepsilon},
\end{equation}
where $\varepsilon>0$ is sufficiently small.

Since 
$$\widetilde{\omega}(y,v)=D^2 v-\widetilde{A}(y,v,Dv)=J^{-1  }  \omega(x,u) (J^{-1})^T,$$
$$\widetilde{B}(y,v,Dv)= J^{-1} B(x,u,Du) (J^{-1})^T,$$
from \eqref{5.3} we have 
\begin{equation*}
\frac{1}{(1+\varepsilon) }  \lambda_{u} \leqslant \lambda_{v} \leqslant(1+\varepsilon) \lambda_{u},
\end{equation*}
\begin{equation*}
\frac{1}{(1+\varepsilon) }  \mu(B) \leqslant \mu(\widetilde{B}). \leqslant(1+\varepsilon) \mu(B).
\end{equation*}
Suppose $u(x)$ is a strictly $\left(\delta, \widetilde{\gamma}_{k}\right)$-admissible solution, i.e.
\begin{equation*}
\mu(B)\le \delta \lambda_u,
\end{equation*}
\begin{equation}\label{5.7}
\lambda_{\min }(\omega(x,u))\ge \widetilde{\gamma}_k \lambda_{\max} (\omega(x,u)),\ x\in \overline{\Omega}.
\end{equation}
From \eqref{5.2}-\eqref{5.7} we obtain
$$  \mu(\widetilde{B})\le (1+\varepsilon)^2 \delta\lambda_v,$$
$$\lambda_{\min} ( \widetilde{\omega}(y,v) )\ge \frac{\widetilde{\gamma}_k}{(1+\varepsilon)^2} \lambda_{\max} ( \widetilde{\omega}(y,v) ).$$
So, $v(y)$ is a strictly $\left(\widetilde{\delta},\widetilde{\widetilde{\gamma}}_k \right)$-admissible solutions to \eqref{1.41}, where 
\begin{equation*}
\widetilde{\delta}= (1+\varepsilon)^2 \delta, \widetilde{\widetilde{\gamma}}_k=\frac{ \widetilde{\gamma}_k }{(1+\varepsilon)^2}.
\end{equation*}
Since $  0<\gamma_k<  \widetilde{\gamma}_k<1,$ where $\gamma_k$ is defined in Definition \ref{dn1.3} and $0<\delta<\delta_{k}<1,$ $\delta_{k}$ is determined in Theorem \ref{dl1.5}, we can assume that $\varepsilon$ is chosen sufficiently small so that
\begin{equation*}
0<\gamma_k< \widetilde{\widetilde{\gamma}}_k <1,\quad  0< \widetilde{\delta} <\delta_k<1.
\end{equation*}

We prove now that the function $\widetilde{F}_k(\widetilde{R}),$ defined by \eqref{1.48}, is $\widetilde{d}$-concave in the sense of \eqref{1.21}.  We rewrite \eqref{1.48} as follows
\begin{equation}\label{5.10}
\widetilde{F}_k\left(y', \widetilde{R}\right)=\log \left[S_{k}(\widetilde{R})+H_{k}\left(y', \widetilde{R}\right)\right], 
\end{equation}
where $H_{k}(y', \widetilde{R})$ is homogeneous of degree $k$ with respect to $\widetilde{R}=\left[\widetilde{R}_{i j}\right]_{n \times n}.$

Suppose $\widetilde{R}=\widetilde{\omega}+\widetilde{\beta} \in D_{\widetilde{\sigma}, \mu(\widetilde{B}),\widetilde{\widetilde{\gamma}}_k}.$
Then we have
$$\widetilde{R}=C^{-1} \widetilde{D} C+\widetilde{\beta}=C^{-1}\left(\widetilde{D}+C \widetilde{\beta} \dot{C}^{-1}\right) C=
C^{-1}(\widetilde{D}+\widetilde{\sigma}) C,$$
where $\widetilde{D}+\widetilde{\sigma}\in D_{\widetilde{\delta}, \mu(\widetilde{B})}.$  So we can assume that $\widetilde{R}=\widetilde{D}+\widetilde{\sigma}, $ $\widetilde{D}=\operatorname{diag} (\widetilde{\lambda}_{1}, \widetilde{\lambda}_{2}, \ldots, \widetilde{\lambda}_{n})>0,$ 
$\widetilde{\lambda}_{\min } \geqslant \widetilde{\widetilde{\gamma}}_{k} \widetilde{\lambda}_{\max },\ \|\widetilde{\sigma}\| \leqslant \mu(\widetilde{B}) \leqslant  \widetilde{\delta}   \widetilde{\lambda}_{\min }.$

From \eqref{5.10} we have 
$$\frac{\partial \widetilde{F}_k(y', \widetilde{R} )}{\partial \widetilde{R} _{ij}} =
\frac{1}{(S_k(\widetilde{R} )+H_k (y', \widetilde{R} )) }.
\frac{\partial (S_k(\widetilde{R} ) + H_k (y', \widetilde{R} )) }{\partial \widetilde{R} _{ij} },$$
\begin{align*}
\frac{\partial^2 \widetilde{F}_k(y', \widetilde{R} )}{\partial \widetilde{R} _{ij}  \partial \widetilde{R} _{\ell m}}=&-\frac{1}{\left( S_k(\widetilde{R} )+H_k (y', \widetilde{R} )   \right)^2}. \frac{ \partial (S_k+H_k)   }{\partial\widetilde{R} _{ij} }.
\frac{ \partial (S_k+H_k)   }{\partial\widetilde{R} _{\ell m} }	\\
&+ \frac{1}{\left( S_k(\widetilde{R} )+H_k (y', \widetilde{R} )   \right)}. \frac{ \partial^2 (S_k+H_k)   }{\partial\widetilde{R} _{ij} \partial\widetilde{R} _{\ell m}}. 
\end{align*}
Then, for $\widetilde{M}=\left[\widetilde{M}_{ij}\right]_{n\times n} \in \mathbb{R}^{n \times n}:$
\begin{equation}\label{5.11}
\begin{split}
	d^{2} \widetilde{F}_{k}\left(y', \widetilde{R}, \widetilde{M}\right) =&-\frac{1}{\left(S_{k}+H_{k}\right)^{2}}\left[ d S_{k}(\widetilde{R}, \widetilde{M})+d H_{k}\left(y', \widetilde{R}, \widetilde{M}\right)\right]^{2} \\
	& +\frac{1}{\left(S_{k}+H_{k}\right)}\left[d^{2} S_{k}(\widetilde{R}, \widetilde{M})+d^{2} H_{k}\left(y', \widetilde{R}, \widetilde{M}\right)\right] .
\end{split}
\end{equation}
We have the following relations:
\begin{equation*}
\frac{1}{S_{k}+H_{k}}=\frac{1}{S_{k}}-\frac{H_{k}}{S_{k}\left(S_{k}+H_{k}\right)},
\end{equation*}
\begin{equation}\label{5.13}
\frac{1}{\left(S_{k}+H_{k}\right)^{2}}=\frac{1}{S_{k}^{2}}-2 \frac{H_{k}}{S_{k}^{2}\left(S_{k}+H_{k}\right)}+\frac{H_{k}^2}{S_{k}^{2}\left(S_{k}+H_{k}\right)^{2}},
\end{equation}
\begin{equation}\label{5.14}
S_{k}(\widetilde{R}) \geqslant S_{k}(\widetilde{D}) \geqslant \binom{n}{k}  (\widetilde{\lambda}_{\min })^{k} \geqslant\binom{n}{k}   (\widetilde{\widetilde{\gamma}}_{k})^k (\widetilde{\lambda}_{\max})^{k}.
\end{equation}
Since the function $H_k(y,\widetilde{R} )$ is a linear combination of $\det \left(  	\widetilde{R}^{(k)}_{i_{1}\cdots i_{k},j_{1}\cdots j_{k}}   \right)$ with coefficients, that are polymomials with respect to $Dh(y')$ of degree at the least 1 and at the most 2, $Dh(y')$ is small, and
 $$|\widetilde{R}_{ij}| \le \delta_{ij} \widetilde{\lambda}_{\max } + | \widetilde{\sigma}_{ij}|\le \delta_{ij} \widetilde{\lambda}_{\max }+\sqrt{n} \widetilde{\delta} \widetilde{\lambda}_{\min }\le (1+ \sqrt{n} \widetilde{\delta}  )\widetilde{\lambda}_{\max },$$
 we have 
 \begin{equation}\label{5.14a}
 	\left| H_k(y', \widetilde{R})\right| \le C_7 |Dh(y')|\left(\widetilde{\lambda}_{\max }\right)^k,
 \end{equation}
$C_7>0$ and does not depend on $y'.$

From \eqref{5.14} and \eqref{5.14a} we can assume that the neighborhood $\widetilde{\Omega}_{\rho}$ is chosen small so that for any $ |y'|\le \rho$
\begin{equation}\label{5.15}
S_{k}(\widetilde{R})+H_{k}\left(y', \widetilde{R}\right) \geqslant C_{8}\left(\widetilde{\lambda}_{\max }\right)^{k}, 
\end{equation}
$C_8>0$ and does not depend on $y'.$

It follows from \eqref{5.11}-\eqref{5.13} that
\begin{equation}\label{5.16}
		\begin{split}		 
d^{2} \widetilde{F}_{k}\left(y', \widetilde{R}, \widetilde{M}\right)=&d^{2}\left(\log S_{k}\right)(\widetilde{R}, \widetilde{M})\\
&+\sum_{i,j,\ell,m=1}^n\left[ \sqrt{ \widetilde{\lambda}_i  \widetilde{\lambda}_j  \widetilde{\lambda}_{\ell}  \widetilde{\lambda}_m}
g_{ij,\ell m} (y', \widetilde{R})\right] \left(\frac{ \widetilde{M}_{ij} }{  \sqrt{ \widetilde{\lambda}_i  \widetilde{\lambda}_j}}\right)
\left(\frac{ \widetilde{M}_{\ell m} }{  \sqrt{ \widetilde{\lambda}_{\ell} \widetilde{\lambda}_m}}\right),
	\end{split}
\end{equation}
where $	g_{ij,\ell m} (y', \widetilde{R})$ are homogeneous of degree $(-2)$ with respect to $\widetilde{R}.$ From \eqref{5.11}-\eqref{5.15} we can assume that for any $i,j,\ell,m$ 
\begin{equation*}
\sup_{\widetilde{\lambda}_{\min } \ge \widetilde{\widetilde{\gamma}}_k \widetilde{\lambda}_{\max }\atop \| \widetilde{\sigma}\|\le \widetilde{\delta}\widetilde{\lambda}_{\min }} \left|  \sqrt{ \widetilde{\lambda}_i  \widetilde{\lambda}_j  \widetilde{\lambda}_{\ell}  \widetilde{\lambda}_m}
g_{ij,\ell m} (y', \widetilde{D}+\widetilde{\sigma})  \right|\le C_9 |Dh(y')|,
\end{equation*}
$C_9>0$ and does not depend on $y', i, j, \ell, m.$

We know from Theorem \ref{dl1.5} that for the function $F_{k}( \widetilde{R})=\log(S_k(\widetilde{R}))$  
when $\widetilde{R}\in D_{\widetilde{\delta},\mu(\widetilde{B}),\widetilde{\widetilde{\gamma}}_k},$ 
where $0<\gamma_{k}<\widetilde{\widetilde{\gamma}}_k<1,$ $0< \widetilde{\delta} <\delta_k<1,$ the estimates \eqref{1.19}, \eqref{1.20} hold, i.e.

\begin{equation*}
d^{2} F_{k}(\widetilde{R},\widetilde{P}) \le -C_1 \frac{| \widetilde{P} |^2 }{\widetilde{\lambda}_{\max }^2 }, \widetilde{P}^T=\widetilde{P}, 
\end{equation*}
\begin{equation}\label{5.19}
d^{2} F_{k}(\widetilde{R}, \widetilde{P}+\widetilde{Q}) \le C_2 \frac{| \widetilde{Q} |^2 }{\widetilde{\lambda}_{\min }^2 }, \widetilde{Q}^T=-\widetilde{Q}. 
\end{equation}
From \eqref{5.16}-\eqref{5.19} it follows that if we choose $\widetilde{\Omega}_{\rho}$ sufficiently small, then we have the following estimates for any $|y'|\le \rho$ 
\begin{equation}\label{5.20}
d^{2} \widetilde{F}_{k}\left(y', \widetilde{R}, \widetilde{P}\right) \leqslant -C_{10} \frac{|\widetilde{P}|^{2}}{\lambda_{\max}^{2}(\widetilde{\omega})}, \widetilde{P}^T  =\widetilde{P},
\end{equation}

\begin{equation}\label{5.21}
d^{2} \widetilde{F}_{k}\left(y', \widetilde{R}, \widetilde{P}+\widetilde{Q}\right) \leqslant C_{11}  \frac{|\widetilde{Q}|^{2}}{\lambda_{\min}^{2}(\widetilde{\omega})}, \widetilde{Q}^T=-\widetilde{Q}^T,
\end{equation}
where $C_{10}>0, C_{11}>0$ depend on $C_1, C_2, C_9, \rho, D h$ and do not depend on $y'$ and $\mu(\widetilde{B}).$
From \eqref{5.21} it is easy to obtain the  following version of $\widetilde{d}$-concavity for the function $\widetilde{F}_{k}(y', \widetilde{R})$ on the set $D_{\widetilde{\delta},\mu(\widetilde{B}), \widetilde{\widetilde{\gamma}}_k}:$
\begin{equation}\label{5.22}
	\begin{split}
	\widetilde{F}_{k}\left(y', \widetilde{R}^{(1)} \right)- \widetilde{F}_{k}\left(y', \widetilde{R}^{(0)} \right) \le &\sum_{ i,j=1}^{n} \frac{\partial \widetilde{F}_{k}\left(y', \widetilde{R}^{(0)} \right)   }{\partial \widetilde{R}_{ij} }(\widetilde{R}_{ij}^{(1)} - \widetilde{R}_{ij}^{(0)})\\
	&+ C_{11}\frac{ |\widetilde{\beta}^{(1)}-\widetilde{\beta}^{(0)}|^2 }{{\lambda}_{\min }^2 (\widetilde{\omega}^{(\tau)} )},\ 0< \tau<1,
	\end{split}
\end{equation}
where $|y'| \le \rho,$ $ \widetilde{R}^{(0)} =\widetilde{\omega}^{(0)} +\widetilde{\beta}^{(0)},$  $ \widetilde{R}^{(1)} =\widetilde{\omega}^{(1)} +\widetilde{\beta}^{(1)}\in D_{\widetilde{\delta},\mu(\widetilde{B}), \widetilde{\widetilde{\gamma}}_k},$ 
$\widetilde{\omega}^{(\tau)}=(1-\tau) \widetilde{\omega}^{(0)}+\tau \widetilde{\omega}^{(1)}.$

\hfill $\square$

\subsection{The H\"older continuity of $D^2 v(y)$}
Since $\|u(x)\|_{C^{2}(\overline{\Omega}) }\leqslant C_{3}$ and $v(y)=u(x)-\varphi(x),$
where $y$ and $x$ are related by \eqref{1.32}, $h(y') \in C^{4},$
then we have
\begin{equation}\label{5.20A}
\|v(y)\|_{C^{2}\left(\overline{\widetilde{\Omega}}_{\rho}\right)} \leqslant C_{12}.
\end{equation}
From \eqref{1.27a} it follows  that
\begin{equation}\label{5.21A}
	\|J^{-1}(y')\|_{C^2(|y'|\le \rho)} \le C'_{12}.
\end{equation}
From \eqref{5.20A}, \eqref{5.21A} and \eqref{1.39} we obtain
\begin{equation}\label{5.22A}
	\|\widetilde{A} (y, v(y), D v(y))\|_{C^{2}\left(\overline{\widetilde{\Omega}}_{\rho}\right)}, 
		\|\widetilde{B} (y, v(y), D v(y))\|_{C^{2}\left(\overline{\widetilde{\Omega}}_{\rho}\right)}, 
			\|\widetilde{f} (y, v(y), D v(y))\|_{C^{2}\left(\overline{\widetilde{\Omega}}_{\rho}\right)} \le C^{"}_{12},
\end{equation}
where $C'_{12}, C^{"}_{12}$ are uniformly bounded when $0<\rho \le \rho_0.$

The matrix $\frac{1}{2}\left(\frac{\partial {F}_k(R)}{\partial R_{ij}}+\frac{\partial F_{k}(R)}{\partial {R_{ji}}}\right)$ satisfies the
ellipticity conditions \eqref{3.6}. But, by definition
$\widetilde{F}_k(y',\widetilde{R})=F_{k}(R)=F_{k}\left(J^{-1} \widetilde{R}\left(J^{-1}\right)^{T}\right), $ where 
$\frac{1}{\sqrt{1+\varepsilon}} \leqslant\left\|J^{-1}\right\| =\left\|(J^{-1})^T\right\| \leqslant \sqrt{1+\varepsilon},$ it follows from \eqref{3.6} that for any $|y'|\le \rho$
$$C_{13}|\xi|^2\le \frac{1}{2} \sum_{i,j=1}^{n} \left(\frac{\partial \widetilde{F}_k(y',\widetilde{R} )}{\partial \widetilde{R}_{ij}}+\frac{\partial \widetilde{F}_{k}(y',\widetilde{R})}{\partial \widetilde{R}_{ji}}\right)\xi_i\xi_j
\le C_{14}|\xi|^2,$$
where $C_{13}, C_{14}$ depend on $n, k, \delta,\widetilde{\gamma}_k, f_0, f_1, \varepsilon$ and do not depend on $y'.$ We have just proved above the strict concavity \eqref{5.20} of the 
function $\widetilde{F}_{k}\left(y', \widetilde{\omega}+\widetilde{\beta}\right)$ when $\widetilde{\beta}$ fixed and the $\widetilde{d}$-concavity \eqref{5.21} of $\widetilde{F}_{k}\left(y', \widetilde{R}\right)$ on the set $D_{\widetilde{\delta}, \mu{(\widetilde{B})}, \widetilde{\widetilde{\gamma}}_k},$ $|y'|\le \rho.$
From the facts listed above with the aid of \eqref{5.20A}-\eqref{5.22A} and \eqref{2.2}-\eqref{2.7N}, applied for $\widetilde{R}(y,v)$ and $\widetilde{F}_k(\widetilde{R}),$ one can prove (\cite{10}, Section 17.8) that 
from the equation $\widetilde{F}_k(y',\widetilde{R}(y,v))=\log \tilde{f}(y,v,Dv)$ in $\widetilde{\Omega}_\rho, v(y)|_{y_n=0}=0$ it follows 
$$
\left\|D^{2} v(y)\right\|_{C^{\alpha}\left(\overline{\widetilde{\Omega}_{\rho}}\right)} \leqslant C_{15}, 
$$
where $C_{15}>0, 0<\alpha<1$ do not depend on $\mu(\widetilde{B}).$ From the last inequality and \eqref{1.32}, \eqref{1.33}, we have: 
$$\| D^{2} u\|_{C^{\alpha}\left(\overline{\Omega}_{\rho}\right)}   \leqslant C'_{15}.$$
\subsection{Proof of Theorem \ref{dl1.11}}
In Sections 3 and 5 we have obtained the following estimates for a strictly $\left(\delta, \widetilde{\gamma}_{k}\right)$-admissible solution $u(x)$ to the problem \eqref{1.1}-\eqref{1.2}:
$$\|u\|_{C^{2}(\overline{\Omega})} \leqslant C_{3},$$
$$\left\|D^{2} u\right\|_{C^{\alpha}\left(\overline{\Omega'}\right)} \leqslant C_{5},\ \Omega' \subset\subset \Omega,$$
$$\left\|D^{2} u\right\|_{C^{\alpha}\left(\overline{\Omega}_{\rho}\right) }  \leqslant C'_{15},\ \Omega_{\rho}= B_{\rho}(x) \cap \Omega,\ x \in \partial \Omega$$ 
From these estimates it follows the desired
inequality \eqref{1.50}:
$$
\|u\|_{C^{2, \alpha}(\overline{\Omega})} \leqslant C_{4},
$$
where $0<\alpha<1,$ $C_4>0$ depend on $n,$ $k,$ $\delta,$ $\widetilde{\gamma}_k,$ $\beta_1,$ $\Omega,$ $A(x,z,p),$ $f(x,z,p),$ $\underline{u}(x),$ $\varphi,$ $\mu_2(B).$
\hfill $\square$

\section{The solvability of the Dirichlet problem}
\subsection{A necessary condition for the existence of a	strictly	$\delta$-admissible solution}
We give here proof for Theorem \ref{dl1.12}.	
Suppose there exists a strictly $\delta$-admissible solution $u(x)$ to the equation \eqref{1.1} and it satisfies 
the conditions (i),(ii) and (iii) of the theorem.

Since $S_k(\omega(x,u)-B(x,u,Du))\ge S_k (\omega(x,u)),$
from \eqref{1.1} it follows that 
$S_k (\omega(x,u)) \le f (x,u,Du).$
But 	 $\omega(x,u)>0,$ $S_k(\omega(x,u)) \ge \binom{n}{k} \lambda_{\min}^k (\omega(x,u))$ 
and 
$\mu(B)\le \delta \lambda_u\le \delta \lambda_{\min} 	(\omega(x,u)),$ then we have
$$\mu(B)\le \delta \left[  \frac{f(x,u,Du) }{\binom{n}{k}}   \right]^{\frac{1}{k}}$$
and consequently 
$$\mu(B)\le \delta \left[  \frac{f_1 }{\binom{n}{k}}   \right]^{\frac{1}{k}},$$	
where $f_1$ is defined by \eqref{1.22a}. The inequality \eqref{1.53} is proved.
\hfill $\square$
\subsection{Some sufficient conditions for unique existence of the strictly  $(\delta,\widetilde{\gamma}_k)$-admissible solution}
We prove here Theorem \ref{dl1.13} on the 
unique solvability of the problem \eqref{1.1}-\eqref{1.2} in the class of strictly $(\delta, \widetilde{\gamma}_k)$-admissible solutions that belong to  $C^{2, \alpha}(\overline{\Omega})$ for some $0<\alpha<1.$ The uniqueness
follows from the comparison principle, Theorem \ref{dl2.2}.

Suppose $B(x, z, p)$ satisfies \eqref{1.55}, \eqref{1.56}.    By using the method of continuity (\cite{10}, Section 17.2) we
will prove the existence of strictly $\left(\delta, \widetilde{\gamma}_{k}\right)$-admissible solution $u(x)$ to the problem \eqref{1.1}-\eqref{1.2}.

Since $S_{k}(R(x, \underline{u})) \geqslant S_{k}(\omega(x, \underline{u})),$ it follows from
the conditions (iii) and (iv) that the function
$\underline{u}(x)$ is also strictly $\left(\delta, \widetilde{\gamma}_{k}\right)$-admissible subsolution to the problem \eqref{1.1}-\eqref{1.2}. Now
for each $t \in[0,1]$ we consider the following Dirichlet problem:

\begin{equation}\label{6.1}
S_k\left[ D^2 u^{(t)} -A(x,u^{(t)}, Du^{(t)} ) -B(x,u^{(t)}, Du^{(t)} )  \right]=f^{(t)}(x,u^{(t)}, Du^{(t)} ) \text{ in } \Omega,
\end{equation}
\begin{equation}\label{6.2}
u^{(t)}=\varphi \text{ on } \partial \Omega,
\end{equation}
where 
\begin{equation}\label{6.3}
f^{(t)}(x,z,p)	=f(x,z,p) e^{(1-t)G(\underline{u})(x)},
\end{equation}	
\begin{equation}\label{6.4}
G[w]		(x) =\log \left(S_k(R(x,w)) \right)-\log f(x,w,Dw).
\end{equation}
From \eqref{6.1}-\eqref{6.4} it follows that the function $u^{(0)}=\underline{u}(x)$ is the solution to the problem \eqref{6.1}-\eqref{6.2} with $t=0$ and if the function $u^{(1)}(x)$ is solution to the problem \eqref{6.1}-\eqref{6.2} when $t=1,$ then $u(x)=u^{(1)}(x)$ is a solution to the problem \eqref{1.1}-\eqref{1.2}.

To study the problem \eqref{6.1}-\eqref{6.2}, for $\varepsilon>0$ we introduce a class $\mathcal{U}^{(\varepsilon)}=\mathcal{U}(k,n,\delta, \widetilde{\gamma}_{k},\varepsilon,\beta_1,B)$ that consists of 
functions $u(x) \in C^{2}(\overline{\Omega}),$ that satisfy the
following conditions:
\begin{equation}\label{6.5}
\lambda_{u}>0,
\end{equation}
\begin{equation}\label{6.6}
\mu(B)<\delta \lambda_{u},
\end{equation}
\begin{equation}\label{6.7}
\mu\left(D_{z} B\right)<\beta_{1} \lambda_{u},
\end{equation}
\begin{equation}\label{6.8}
\gamma_{u}>\widetilde{\gamma}_{k}+\varepsilon,
\end{equation}
where $\lambda_u$ and $\gamma_{u}$ are defined by \eqref{1.11} and \eqref{1.22} respectively.

We note that if $\varepsilon_1>\varepsilon_2>0$ then $\mathcal{U}^{(\varepsilon_1)}\subset  \mathcal{U}^{(\varepsilon_2)}.$
It is obvious that $\mathcal{U}^{(\varepsilon)}$ is open in $C^{2}(\overline{\Omega}).$ From \eqref{1.53a}-\eqref{1.56}  it follows that $\underline{u}(x) \in \mathcal{U}^{(\varepsilon_0)}.$ If $u^{(t)} \in \mathcal{U}^{(\varepsilon)}$ and it is a solution to the problem \eqref{6.1}-\eqref{6.2}, then from \eqref{6.5}-\eqref{6.8} we see that it is also a strictly  $\left(\delta, \widetilde{\gamma}_{k}\right)$-admissible solution.
\begin{lemma}\label{bd6.1}
Assume that all conditions of Theorem \ref{dl1.13} are fulfilled. Then, $A(x, z, p),$
$B(x, z, p),$
$f^{(t)}(x, z, p)$ satisfy all assumptions of Proposition \ref{md1.7} and of Theorems \ref{dl1.8}, \ref{dl1.11}, and the function $\underline{u}(x)$ is a strictly $\left(\delta, \widetilde{\gamma}_{k}\right)$-admissible subsolution to the all problems   \eqref{6.1}- \eqref{6.2}.	
\end{lemma}
\begin{proof}
Since $f^{(t)}(x,z,p)	=f(x,z,p) e^{(1-t)G(\underline{u})(x)},$ 
then $f^{(t)}(x,z,p)>0$ and 
$$\inf_{\mathcal{D}}\left[ \frac{ D_z f^{(t)}(x,z,p)}{f^{(t)}(x,z,p)}  \right] =\inf_{\mathcal{D}} \left[ \frac{ D_z f (x,z,p)}{f (x,z,p)}  \right] \ge \frac{k \delta}{(1+\delta^2)} \beta_{1}. $$	
We show that the function $\underline{u}(x)$ is strictly $\left(\delta, \widetilde{\gamma}_{k}\right)$-admissible subsolution to all equations \eqref{6.1}. 

Indeed, since $G\left[\underline{u}\right](x)\ge 0,$ we have for $0\le t \le 1:$
\begin{align*}
	S_{k}(R(x, \underline{u}))&=f(x, \underline{u}, D \underline{u}) \cdot \frac{S_{k}(R(x, \underline{u}))}{f(x, \underline{u}, D \underline{u})} \\
	&=	f(x, \underline{u}, D \underline{u})  e^{G\left[ \underline{u}\right](x)}\ge f(x, \underline{u}, D \underline{u})   e^{(1-t)G\left[ \underline{u}\right](x)}=f^{(t)}(x, \underline{u}, D \underline{u}).
\end{align*}
Moreover, from \eqref{1.55}, \eqref{1.56}, \eqref{6.6}, \eqref{6.7} it follows that
$$\mu(B)<\delta \min \left(\lambda_{u}, \lambda_{\underline{u}}\right),$$
$$\mu\left(D_{z} B\right)<\beta_{1} \min \left(\lambda_{u}, \lambda_{\underline{u}}\right).$$
\end{proof}
From Lemma \ref{bd6.1}, Proposition \ref{md1.7}, Theorems \ref{dl1.8} and \ref{dl1.11} we obtain 

\begin{corollary}\label{hq6.2}
Suppose $u^{(t)}$ is a strictly $(\delta,\widetilde{\gamma}_k)-$admissible solution to the problem \eqref{6.1}-\eqref{6.2}. Then there
exist $M_{0}>0, M_{1}>0,0<\alpha<1, C_{4}>0$  that 
depend on $n, k,$ $ \delta,$  $\widetilde{\gamma}_{k},$ $\beta_{1},$ $\Omega,$ $A(x,z, p),$ $f(x, z, p),$ $\underline{u}(x),$ $\varphi,$ $\mu_2(B),$ and do not depend on $t$ such that

$$\sup_{\overline{\Omega}}\left|u^{(t)}\right| \leqslant M_{0}, \quad 
\sup_{\overline{\Omega}}\left|D u^{(t)}\right| \leqslant M_{1},$$
\begin{equation}\label{6.9}
	\| u^{(t)}\|_{C^{2,\alpha}(\overline{\Omega})} \le C_4,
\end{equation}
\begin{equation}\label{6.10}
	\lambda_{u^{(t)}} \ge \widetilde{\gamma}_k\left[ \frac{(1+\delta^2)^{-\left[\frac{k}{2}\right]}}{\binom{n}{k}} f_0 \right]^{\frac{1}{k}}.
\end{equation}
\end{corollary}

Here, to get \eqref{6.10} we have used from \eqref{1.22} the fact that
$$f_0= \inf_{x \in \overline{\Omega} \atop |z|\le M_0, |p|\le M_1}  f(x,z,p)
\le \inf_{x \in \overline{\Omega} \atop |z|\le M_0, |p|\le M_1}  f(x,z,p) e^{(1-t)G\left[ \underline{u}\right](x)}.$$

We rewrite the problem \eqref{6.1}-\eqref{6.2}  as follows
\begin{equation}\label{6.11}
	\log \left(S_{k}\left(R\left(x, u^{(t)}\right)\right)\right)-\log f\left(x, u^{(t)}, D u^{(t)}\right)  =(1-t) G[\underline{u}](x) \ \text { in } \Omega, 
\end{equation}
\begin{equation}\label{6.12}
	u^{(t)}=\varphi \ \text { on } \partial \Omega. 
\end{equation}
We consider the operator:
$$G[u](x): C^{2,\alpha}\left(\overline{\Omega}\right) \to C^{0,\alpha}\left(\overline{\Omega}\right),$$
where $G[u](x)$ is defined by \eqref{6.4}, which is connected to the left-hand side of \eqref{6.11} and $0<\alpha<1$ is the same as in \eqref{6.9}.

\begin{lemma}\label{bd6.3}
Suppose 	$u \in C^{2,\alpha}\left(\overline{\Omega}\right)$ is a strictly $(\delta,\widetilde{\gamma}_k)-$admissible solution to the 	problem \eqref{6.11}-\eqref{6.12}. Then the operator $	G[u](x)$ is Frechet continuously differentiable at 	$u$ and its differential $G_u$ is defined as follows 
$$G_u: C^{2,\alpha}_0\left(\overline{\Omega}\right)  \to C^{0,\alpha}\left(\overline{\Omega}\right),$$
where $C^{2,\alpha}_0\left(\overline{\Omega}\right) =\left\{  h\in C^{2,\alpha}\left(\overline{\Omega}\right); h=0 \text{ on } \partial \Omega \right\}.$
$$G_u(h)= \sum_{ i,j=1}^n a^{ij}(x) D_{ij}h+ \sum_{ i=1}^n b^i (x) D_ih+c(x) h,$$
$$a^{ij}(x)=\frac{1}{2} \left[ F^{ij}[u](x)+ F^{ji}[u](x)\right],\ i,j=1,\cdots,n,$$
$$b^i(x)=-\sum_{\ell, m=1}^n F^{\ell m}[u](x) D_{p_i} (A_{\ell m}+B_{\ell m})(x,u,Du)-\left( \frac{D_{p_i}f}{f}\right)(x,u,Du),\ i=1,\cdots,n,$$
$$c(x)=-\sum_{\ell, m=1}^n F^{\ell m}[u](x) D_{z} (A_{\ell m}+B_{\ell m})(x,u,Du)-\left( \frac{D_{z}f}{f}\right)(x,u,Du),  $$
$$F^{ij}[u](x)=\frac{\partial F_k(R(x,u))}{\partial R_{ij}}, F_k(R)=\log(S_k(R)).$$
The operator $G_{u}(h)$ is uniformly elliptic on $\overline{\Omega},$ all the coefficients $a^{ij}, b^i, c$ are from $C^{0,\alpha}\left(\overline{\Omega}\right)$ and $c(x)\le 0.$ Moreover, it is invertible.
\end{lemma}
\begin{proof}
Since $u(x)\in C^{2,\alpha}(\overline{\Omega})$ and it is a strictly $\left(\delta, \widetilde{\gamma}_{k}\right)$-admissible solution to the problem \eqref{6.1}-\eqref{6.2}, 
then the uniform ellipticity of the operator $G_{u}(h)$ follows from \eqref{3.6}.  Due to $A(x,z,p),$ $ B(x,z,p),$ $ f(x,z,p)\in C^3(\mathcal{D}),$ $u(x)\in C^{2,\alpha}\left(\overline{\Omega}\right),$ then the coefficients $a^{ij}(x), b^i(x), c(x)$ are from $C^{0,\alpha}\left(\overline{\Omega}\right).$ As in the proof of the comparison principle, from the assumptions on $D_z A,$ $ D_z B$ and $D_z f,$ the assertion $c(x)\le 0$ can be verified, from which it follows that $G_u$ is invertible.
\end{proof}
We rewrite the problem \eqref{6.11}-\eqref{6.12} in the form
\begin{equation}\label{6.12a}
H(u^{(t)}, t)=0 \ \text{ in }\Omega,\quad  u^{(t)}=\varphi \text{ on } \partial \Omega
\end{equation}
where $	H: C^{2,\alpha}\left(\overline{\Omega}\right)\times [0,1]\to C^{0,\alpha}\left(\overline{\Omega}\right), $
\begin{equation}\label{6.14}
	H(u^{(t)},t)= G[u^{(t)}](x)-(1-t)G[\underline{u}](x).
\end{equation}
We consider a set of solutions to the problem \eqref{6.12a}  as follows
$$\mathcal{V}^{(\varepsilon)}= \mathcal{U^{(\varepsilon)}}\cap C^{2,\alpha}\left(\overline{\Omega}\right), $$
where $0<\alpha<1$ as in \eqref{6.9}, fixed and is the same for all $u^{(t)},$ $0\le t\le 1.$

We introduce the following set 
\begin{equation*}
I=\left\{t\in [0,1]: \exists u^{(t)}\in \mathcal{V}^{(\varepsilon)}, \varepsilon= \varepsilon(u^{(t)})>0, H(u^{(t)},t)   =0, u^{(t)}=\varphi \text{ on } \partial \Omega  \right\}.
\end{equation*}
The solvability of the problem \eqref{6.11}-\eqref{6.12}
is equivalent to the fact that $t \in I.$
When $t=0$ the function $u^{(0)} =\underline{u}$ is a
solution to \eqref{6.12a}, i.e.
$$
H(u^{(0)}, 0)=0.
$$
This means that $t=0\in I$ and $I \neq \emptyset.$ The following lemma
shows that $I$ is "open".
\begin{lemma}\label{bd6.4N}
Suppose $t'\in I,$ 
\begin{equation}\label{6.16}
	u^{(t')}		\in \mathcal{V}^{(\varepsilon')},
\end{equation}
and $\varepsilon'>\varepsilon">0.$ Then there exists $\tau'>0$ such that $[t',t'+\tau'] \subset I$ and 
\begin{equation}\label{6.17}
	u^{(t)}		\in \mathcal{V}^{(\varepsilon")}		
\end{equation}
for any $t\in [t',t'+\tau'].$ Moreover, all $u^{(t)}$ are in some $C^{2,\alpha}(\overline{\Omega})$-neighborhood of $u^{(t')}$ and $u^{(t)}$ is continuous mapping from $[t',t'+\tau']$ to $C^{2,\alpha}(\overline{\Omega}).$ 
\end{lemma}
\begin{proof}
From \eqref{6.14} and Lemma \ref{bd6.3} it follows that the derivative $H_{u^{(t')}	}= G_{u^{(t')}	}$ is invertible, so we can apply the implicit function Theorem to conclude that there exist $\tau'>0$  
and continuous mapping $u^{(t)}$ from $[t',t'+\tau']$ to $C^{2,\alpha}(\overline{\Omega})$ such that
$$H( u^{(t)}  ,t)=0,\ t\in [t',t'+\tau'],\quad  u^{(t)} =\varphi \text{ on } \partial \Omega.$$
We have
 $$\lambda_{\min}(\omega(x,u))=\inf_{|\xi|=1} \sum_{ i,j=1}^n \left[ D_{x_ix_j}u(x)-A_{ij}(x,u,Du)\right] \xi_i \xi_j ,$$ 
$$\lambda_{\max}(\omega(x,u))=\sup_{|\xi|=1}\sum_{ i,j=1}^n \left[ D_{x_ix_j}u(x)-A_{ij}(x,u,Du)\right] \xi_i \xi_j $$
 and
 $A(x,z,p)\in C^3(\mathcal{D}).$  So from \eqref{6.16}, \eqref{1.22} and $\varepsilon'>\varepsilon''>0$ it follows that, if there is a necessity,  
    we may decrease $\tau$ in that way so that \eqref{6.17} is satisfied.
\end{proof}

The desired conclusion of the theorem will be derived from the following lemma.

\begin{lemma}\label{bd6.4}
The assertion 
$$I=[0,1]$$
is true.
\end{lemma}

\begin{proof}
We now apply consecutively Lemma \ref{bd6.4N}. Since $u^{(0)}=\underline{u} \in \mathcal{V}^{\left(\varepsilon_{0}\right)},$ $\varepsilon"= \varepsilon_{1}=\varepsilon_{0}-\frac{\varepsilon_{0}}{4}<\varepsilon_{0}=\varepsilon',$ then   for $t'=0$ there exists $\tau_{1}>0$ such
that if $t_{1}=t'+\tau_{1},$ then $\left[0, t_{1}\right] \subset I$ and
$u^{(t)} \in \mathcal{V}^{\left(\varepsilon_{1}\right)},$ for any $t \in\left[0, t_{1}\right].$
Now we choose $t'=t_{1}$ and $\varepsilon"= \varepsilon_{2}=\varepsilon_{0}-\left(\frac{\varepsilon_{0}}{4}+\frac{\varepsilon_{0}}{8}\right)<\varepsilon_{1}=\varepsilon'.$
Then there exists $\tau_{2}>0$ such that if $t_{2}=t_{1}+\tau_{2}$ then $\left[t_{1}, t_{2}\right] \subset I$ and \eqref{6.17} yields $u^{(t)} \in \mathcal{V}^{\left(\varepsilon_{2}\right)},$ for any $t \in\left[t_{1}, t_{2}\right].$

We set for $m=1,2, \cdots$
$$
\varepsilon"=\varepsilon_{m}=\varepsilon_{0}-\left(\frac{\varepsilon_{0}}{4}+\frac{\varepsilon_{0}}{8}+\cdots+\frac{\varepsilon_{0}}{2^{m+1}}\right), \varepsilon'=\varepsilon_{m-1}.
$$
Then $\varepsilon'> \varepsilon">0$ and we choose $t'=t_{m-1}.$ There exists $\tau_m>0$ such that if $t_m=t_{m-1}+\tau_m$ then $[t_{m-1},t_m]\subset I$ and 	
$u^{(t)} \in \mathcal{V}^{\left(\varepsilon_{m}\right)},$ for any $t \in\left[t_{m-1}, t_{m}\right].$
We can continue this process many times   until $t_{m}<1.$ We set  
$$
t^{*}=\sup_{m \geqslant 1} t_{m}.
$$

We show that $t^*\in I.$
Indeed, we consider the sequence $\left\{u^{ \left(t_{m}\right)}\right\}  .$
Since  $\mathcal{V}^{\left(\varepsilon'\right)} \subset \mathcal{V}^{\left(\varepsilon"\right)}$ if $\varepsilon'>\varepsilon">0$ and $\varepsilon_{m-1}>\varepsilon_{m}>\frac{\varepsilon_0}{2}$
then
\begin{equation}\label{6.18}
	u^{\left(t_{m}\right) }\in \mathcal{V}^{\left(\frac{\varepsilon_{0}}{2}\right)},\ m = 1,2, \cdots 
\end{equation}
From \eqref{6.9}, \eqref{6.10}, \eqref{6.11} we have
\begin{equation}\label{6.16a}
	\left\|u^{\left(t_{m}\right)}\right\|_{ C^{2,\alpha}\left(\overline{\Omega}\right)  } \leqslant C_{4},
\end{equation}
\begin{equation}\label{6.20}
	\lambda_{u^{\left(t_{m}\right)}} \ge \widetilde{\gamma}_k \left[ \frac{ (1+\delta^2)^{-\left[\frac{k}{2}\right]}  }{\binom{n}{k}} f_0 \right]^{\frac{1}{k}},
\end{equation}
and
\begin{equation}\label{6.21}
	G\left[ u^{\left(t_{m}\right)} \right](x)-(1-t_m) G\left[ \underline{u}\right](x)=0.
\end{equation}
From \eqref{6.16a} it follows that there exist $\left\{t_{m'} \right\}\subset \left\{t_{m} \right\}$ and $u(x)\in C^{2,\alpha}\left(\overline{\Omega}\right)  $ such that $t_{m'} \to t^*, $ $u^{(t_{m'} )}\to u(x)$ as $m'\to \infty$ in $C^{2,\alpha}\left(\overline{\Omega}\right).$

Then from  \eqref{6.20}, \eqref{6.21} we obtain
\begin{equation}\label{6.22}
	\lambda_{u} \ge \widetilde{\gamma}_k \left[ \frac{ (1+\delta^2)^{-\left[\frac{k}{2}\right]}  }{\binom{n}{k}} f_0 \right]^{\frac{1}{k}},
\end{equation}
\begin{equation}\label{6.23}
	G\left[ u\right](x)-(1-t^*) G\left[ \underline{u}\right](x)=0.
\end{equation}
But from \eqref{6.18} we have
\begin{equation}\label{6.24}
	\gamma_{u} \geqslant \widetilde{\gamma}_{k}+\frac{\varepsilon_{0}}{2}>\widetilde{\gamma}_{k}+\frac{\varepsilon_{0}}{4}.
\end{equation}
We will verify the conditions \eqref{6.6}, \eqref{6.7}. 

From \eqref{1.55}, \eqref{1.56} and \eqref{6.22} it follows that
\begin{equation}\label{6.25}
	\mu(B)<\delta \lambda_{u},
\end{equation}
\begin{equation}\label{6.26}
	\mu\left(D_{z} B\right)<\beta_{1} \lambda_{u}.
\end{equation}
Therefore, the conditions \eqref{6.6}, \eqref{6.7} are satisfied.   Since $u(x) \in C^{2, \alpha}(\overline{\Omega}),$ from \eqref{6.23}-\eqref{6.26}   
it follows that $u^{\left(t^{*}\right)}=u \in \mathcal{V}^{\left(\frac{\varepsilon_{0}}{4}\right)}$ and
$t^{*} \in I.$

The case $t^*<1$ is impossible, because if $t^{*}<1$ then we can apply again Lemma \ref{bd6.4N} with $t'=t^*, \varepsilon'=\frac{\varepsilon_{0}}{4}, \varepsilon''=\frac{\varepsilon_{0}}{8}<\varepsilon'$ and deduce that there exists $\tau>0$ such that $[t^*,t^*+\tau]\subset I.$ Hence $t^*=1$ and the function
$$u(x)=u^{(1)}\in \mathcal{V}^{\left(\frac{\varepsilon_{0}}{4}\right)}$$
is a strictly $\left(\delta, \widetilde{\gamma}_k+ \frac{\varepsilon_{0}}{4}\right)$-admissible solution to the Dirichlet problem \eqref{1.1}-\eqref{1.2}. The lemma and Theorem \ref{dl1.13} are proved.
\end{proof}
\begin{remark}[On simplified sufficient conditions]\label{nx6.5}
Since $u^{(0)}=\underline{u},$ from \eqref{6.10} it follows that
\begin{equation}\label{6.23a}
	\lambda_{\underline{u}} \ge \widetilde{\gamma}_k \left[ \frac{ (1+\delta^2)^{-\left[\frac{k}{2}\right]}  }{\binom{n}{k}} f_0 \right]^{\frac{1}{k}} =\lambda_{*},
\end{equation}
where $\lambda_{*}$ is defined by \eqref{1.38a}.

Then, from \eqref{6.23a}, \eqref{1.55}-\eqref{1.56} we deduce that for the existence of strictly $\left(\delta, \widetilde{\gamma}_{k}\right)$-admissible solution to the problem \eqref{1.1}-\eqref{1.2},  the matrices $B(x, z, p)$ must satisfy the following simplified sufficient conditions:
\begin{equation}\label{6.27}
	\mu(B)  <\delta \lambda_{*}, 
\end{equation}
\begin{equation}\label{6.28}
	\mu\left(D_{z} B\right)<\beta_{1} \lambda_{*}.
\end{equation}
\end{remark}
The condition \eqref{6.27} is stricter than the necessary condition \eqref{1.53}.
\begin{remark}[On the choice of $\widetilde{\gamma}_k$ and $\delta$] From \eqref{6.23a}, \eqref{6.27} it follows that to have a broader class of the  matrices $B(x,z,p)$ we must increase $\widetilde{\gamma}_k$ and $\delta$ as much as possible. 
	If $k\in \{2,3,n-1,n\},$ then we have to determine $\gamma_{k}$ before we do it for $\widetilde{\gamma}_k.$	
	The parameter $\widetilde{\gamma}_k,$ $0<\gamma_{k} < \widetilde{\gamma}_k<1,$ depends on the choice of the subsolution $\underline{u}(x).$ When $\widetilde{\gamma}_k$ has been chosen, the parameter $\delta_{k},$ $0<\delta_{k}<1,$ is determined as in Theorem \ref{dl1.5}. Then we should choose $\delta=\delta_{k}-\varepsilon_{1},$ where $\varepsilon_{1}>0$ is sufficiently small such that $0<\delta<\delta_k.$	
\end{remark}
\section{An example}
\subsection{A $k$-Hessian type equation in an ellipsoid}
Consider the following problem with $2\le k\le n$ 
\begin{equation}\label{7.1}
S_k\left(D^2 u-A(x,u,Du)-B(x,u,Du)  \right)=f(x,u,Du) \text{ in } \Omega,
\end{equation}
\begin{equation}\label{7.2}
u=0 \text{ on } \partial \Omega,
\end{equation}
where 
\begin{equation}\label{7.2a}
A(x,z,p)=(\arctan z). \frac{|p|^2}{(1+|p|^2)^{\frac{3}{4}}} E_n,	
\end{equation}
$$f(x,z,p)= e^z (1+|p|^2)^m, 0\le m<\frac{k}{2},$$
$$\Omega=\left\{ x\in \mathbb{R}^n: \sum_{j=1}^{n} \mu^2_j x^2_j<1, \mu_j>0 \right\}.$$
We will show how to determine $\gamma_{k}$ for $k\in \{2,3,n-1,n\}$ and how to construct a strictly $\widetilde{\gamma}_k$-admissible subsolution $\underline{u}(x),$ where $0<\gamma_{k}< \widetilde{\gamma}_k<1,$ in this concrete case.

We set $\mu_{\min}=\min_{1\le j\le n}\mu_j, \mu_{\max}= \max_{1\le j\le n} \mu_j,$ $\gamma_{\Omega}=\frac{\mu_{\min}^2 }{\mu_{\max}^2}$ and assume that
\begin{equation}\label{7.4}
\gamma_{k} < \gamma_{\Omega} \le 1		,  \text{ if  } 4\le k \le n-2,
\end{equation}
where $\gamma_{k}, 4\le k \le n-2,$ is defined  in \eqref{1.16} and \eqref{1.17}. In the cases $k \in \{2,3,n-1,n\}$ we may choose $\gamma_{k}=\gamma_{\Omega}-3\varepsilon_{0},$ where $\varepsilon_{0}>0$ is sufficiently small such that $\gamma_{k}>0.$

Then for all $k,$ $2 \le k \le n,$ from \eqref{7.4} we can choose
\begin{equation}\label{7.4a}
\widetilde{\gamma}_k= \gamma_{\Omega}-2\varepsilon_{0},
\end{equation}
where $\varepsilon_{0}>0$ is assumed to be sufficiently small such that $0<\gamma_{k} <\widetilde{\gamma}_k<1.$ For $x\in \overline{\Omega}$ we set 
$$v(x)= \sum_{ j=1}^n \mu^2_j x^2_j-1,$$
\begin{equation}\label{7.4b}
	\underline{u}(x)= \frac{c}{2} v(x),\ c>0.
\end{equation}

We show that if $c> 0$ is chosen sufficiently large, then $\underline{u}(x)$ is a strictly 
$(\widetilde{\gamma}_{k}+\varepsilon_{0})$-admissible subsolution to the equation
\begin{equation}\label{7.6}
S_k\left[D^2 u-A(x,u,Du) \right]=f(x,u,Du)\ \text{ in }\Omega.
\end{equation}
This means that $\underline{u}(x)$ satisfies the following conditions
\begin{equation}\label{7.7}
S_{k}\left[D^{2}\underline{u}  -A(x, \underline{u}, D \underline{u}) \right] \geqslant f\left(x, \underline{u}, D \underline{u} \right) \text { in } \Omega,
\end{equation}
\begin{equation}\label{7.8}
\lambda_{\underline{u}}>0,
\end{equation}
\begin{equation}\label{7.9}
\gamma_{\underline{u}}> \widetilde{\gamma}_k +\varepsilon_{0},\  \varepsilon_{0}>0.
\end{equation}
Indeed, from \eqref{7.4b} we have
$$D\underline{u}=c\left(\mu_{1}^2 x_{1}, \cdots, \mu_{n}^2 x_{n}\right),$$
$$D^{2} \underline{u}=c \operatorname{diag}\left(\mu_{1}^2, \cdots, \mu_{n}^2\right).$$
Since 
\begin{equation*}
A(x,\underline{u},D \underline{u})=\frac{(\arctan \underline{u}(x)) |D\underline{u}(x)|^2}{(1+|D\underline{u}(x)|^2)^{\frac{3}{4}}} E_n
\end{equation*}
and $-\frac{c}{2}\le \underline{u}(x)\le 0$ in $\overline{\Omega},$ $ |D\underline{u}(x)|\le  c \sqrt{n} \mu_{\max},$ then
$$0\le -A(x,\underline{u},D \underline{u})\le \frac{\pi|D\underline{u}(x)|^{2}  }{2(1+|D\underline{u}(x)|^2)^{\frac{3}{4}}}E_n\le \frac{\pi}{2} \sqrt{c\sqrt{n}\mu_{\max}} E_n.$$
Hence, with $\omega(x, u)=D^{2} u-A(x, u, D u)$ we have
\begin{equation}\label{7.11}
\lambda_{\min }(\omega(x, \underline{u})) \geqslant c \mu_{\min }^2,
\end{equation}
\begin{equation}\label{7.12}
\lambda_{\max }(\omega(x, \underline{u})) \leqslant c \mu_{\max}^2+\frac{\pi}{2} \sqrt{c\sqrt{n} \mu_{\max}}.
\end{equation}

From \eqref{7.11}, \eqref{7.12} and \eqref{7.4a} it follows that if we choose $c$ so that $c>c_1,$ where 
\begin{equation*}
c_1=\left(\frac{\pi (\widetilde{\gamma}_k+\varepsilon_{0})}{2\varepsilon_{0}}  \right)	^2 \frac{\sqrt{n}}{\mu_{\max}^3},
\end{equation*}
then 
\begin{equation*}
\gamma_{\underline{u}}=\inf_{x \in \overline{\Omega}}\left(\frac{\lambda_{\min }(\omega(x, \underline{u}))}{\lambda_{\max }\left(\omega\left(x, \underline{u}\right)\right)}\right)>\widetilde{\gamma}_k+  \varepsilon_{0}
\end{equation*}
and \eqref{7.9} is satisfied. From\eqref{7.11} we have 
$\lambda_{\underline{u}}\ge c \mu_{\min}^2$ and \eqref{7.8} holds.

Now we consider the condition \eqref{7.7}. Since $\underline{u}\le 0, -A\ge 0$ and $A(x,z,p)$ is a multiple of $E_n,$ $A(x,\underline{u},D \underline{u} )$ and $D^2 \underline{u}$ commute, we have
\begin{equation*}
\begin{split}
	S_{k}\left[D^{2}\underline{u}  -A(x, \underline{u}, D \underline{u}) \right] &- f\left(x, \underline{u}, D \underline{u} \right)\\
	&\ge S_k(D^2 \underline{u})		-f\left(x, \underline{u}, D \underline{u} \right)\\
	&=c^k\left[ \sigma_k(\mu_1^2,\cdots,\mu_n^2)- \frac{e^{\underline{u}} (1+|D\underline{u}|^2)^m}{c^k} \right]\\
	&\ge c^k\left[ \sigma_k(\mu_1^2,\cdots,\mu_n^2)- \frac{  (1+c^2 n \mu_{\max}^2)^m}{c^k} \right]		.
\end{split}
\end{equation*}
Since $0\le 2m <k,$ the equation 
\begin{equation}\label{7.16}
\frac{  (1+c^2n \mu_{\max}^2)^m}{c^k}=\sigma_k (\mu_1^2,\cdots,\mu_n^2)
\end{equation}
has at least one positive root. We denote by $c_2$ the largest positive root of the equation \eqref{7.16}. Then, when $c>c_2$ we have 
$$S_{k}\left[D^{2}\underline{u}  -A(x, \underline{u}, D \underline{u}) \right] \ge  f\left(x, \underline{u}, D \underline{u} \right)$$
and \eqref{7.7} holds. Hence, if $c>\max(c_1,c_2)$ 
then the function $\underline{u}(x)=c v(x)$ is a strictly  $(\widetilde{\gamma}_k+\varepsilon_{0})$-admissible subsolution to the equation \eqref{7.6}.
The function $\underline{u}(x)$ is also a strictly  $(\widetilde{\gamma}_k+\varepsilon_{0})$-admissible
subsolution to the problem \eqref{7.1}-\eqref{7.2}  for
any skew-symmetric matrix $B(x, z, p) \in B C^{3}(\mathcal{D}).$

Suppose $0<\delta < \delta_{k},$  where $0<\delta_{k}<1$
is determined as in Theorem \ref{dl1.5}. It is
obvious that the matrix $A(x, z, p)$ satisfies the
condition (i) of Theorem \ref{dl1.13}. The 
function $f(x, z, p)>0$ and $\frac{D_{z} f(x, z, p)}{f(x, z, p)}=1.$
So we choose $\beta_{1}=\frac{(1+\delta^{2})}{k \delta}.$

By $A(x, z, p), \underline{u}(x)$ and $\varphi=0$ we determine
$M_{0}>0, M_{1}>0$ as in Theorem \ref{dl1.8}. Then
$$f_0= \inf_{x \in \overline{\Omega} \atop |z|\le M_0, |p|\le M_1} f(x, z, p) =e^{-M_0}$$
and 
$$	\lambda_{*} = \widetilde{\gamma}_k \left[ \frac{ (1+\delta^2)^{-\left[\frac{k}{2}\right]} e^{-M_0} }{\binom{n}{k}}   \right]^{\frac{1}{k}} .$$

Theorem \ref{dl1.13} and Remark \ref{nx6.5} state that the problem \eqref{7.1}-\eqref{7.2} has unique strictly $(\delta,\widetilde{\gamma}_k)$-admissible
solution $u(x),$ that belongs also to $C^{2,\alpha}\left(\overline{\Omega}\right)$ for some $0<\alpha <1,$ if the skew-symmetric matrices $B(x,z,p),$ by \eqref{6.27}, \eqref{6.28}, satisfy the  following conditions:
\begin{equation}\label{7.17}
\mu(B) < \delta \lambda_{*},
\end{equation}
\begin{equation}\label{7.18}
\mu\left(D_{z} B\right)<\frac{(1+\delta^{2})}{k \delta} \lambda_{*}.
\end{equation}
The parameter $0<\alpha<1$ depends on $n, k,$ $ \delta,$  $\widetilde{\gamma}_{k},$   $\mu_2(B).$ 

\begin{remark}
	Since the matrix $A(x,z,p),$ defined by \eqref{7.2a}, does not satisfy the regularity condition (\cite{8}):
	\begin{equation}\label{T7.16}
	\sum_{ i,j, \ell,m=1}^n \frac{\partial A_{ij}(x,z,p)}{\partial p_{\ell} \partial p_m} \xi_i \xi_j \eta_{\ell} \eta_{m} \ge 0,\ (x,z,p)\in \mathcal{D},\ \xi,\eta\in\mathbb{R}^n,\ \xi\perp \eta,
	\end{equation}
	then the equations \eqref{7.1} have not yet been considered in
	\cite{9} of the case $B(x,z,p)=0$ and in \cite{11} of the case $k=n,$ $B(x,z,p)\ne 0.$
So, the result of the Theorem \ref{dl1.13}
 is new even for symmetric $k$-Hessian type equations and nonsymmetric Monge-Ampère type equations. We note that it is the geometric structure condition \eqref{7.4}, that allows one to drop out the condition \eqref{T7.16} for the matrix $A(x,z,p).$
\end{remark}

\subsection{The case $k=2$}
We consider the same equation \eqref{7.1}-\eqref{7.2}, but in the case $k=2.$ It is well-known that if $\widetilde{\omega}, \widetilde{\beta}$ are any matrices of size $2 \times 2$ with $\widetilde{\omega}^{T}=\widetilde{\omega}, \widetilde{\beta}^{T}=-\beta,$ then
$$
\det(\widetilde{\omega}+\widetilde{\beta})=\det\widetilde{\omega}+\det \widetilde{\beta}.
$$
Since $D^{2} u(x)+A(x, u, D u)$ is symmetric, $B(x, u, D u)$ is skew-symmetric, from the assertions (vi), (ix) of Proposition \ref{md4.1}, it follows that the equation \eqref{7.1} becomes the following
\begin{equation}\label{7.19}
S_{2}\left(D^{2} u-A(x, u, D u)\right)=f(x, u, D u)-S_{2}(B(x, u, D u)),\ x\in \Omega,
\end{equation}
where for $B(x,z,p)=\left[B_{ij}(x,z,p)\right]_{n\times n},$ $B^T=-B$ we have 
$$S_2(B(x,z,p))= \sum_{ i<j} B_{ij}^2(x,z,p).$$

That means, we have reduced a nonsymmetric 2-Hessian type equation to a symmetric one with a new right-hand side. Suppose $A(x, z, p)$  and $\Omega$ are the same as in the problem \eqref{7.19},\eqref{7.2} and $\widetilde{\gamma}_2, \gamma_{2}$ are chosen as the same as above, i.e.
$$0<\gamma_{2}=\gamma_{\Omega}-3\varepsilon_{0} < \widetilde{\gamma}_2 = \gamma_{\Omega}-2\varepsilon_{0}<1,\  \varepsilon_{0}>0,$$
where
$\gamma_{\Omega}= \frac{\mu_{\min}^2}{\mu_{\max}^2}.$

 We assume that the function
  $$
g(x, z, p)=f(x, z, p)-S_{2}(B(x, z, p))
$$
satisfies the following conditions:  
\begin{equation}\label{7.15}
g(x, z, p)>0 \text{ in } \mathcal{D},
\end{equation}
\begin{equation}\label{7.21}
D_{z} g(x, z, p) \geqslant 0 \text{ in } \mathcal{D},
\end{equation}
\begin{equation}\label{7.22}
g\left(x, z,p \right) \leqslant C\left(1+|p|^{2}\right)^h,\ 0 \leqslant h<1,\ C>0.
\end{equation}

Then, as for the problem \eqref{7.1}-\eqref{7.2}, we can show that the function $\underline{u}(x)=\frac{c}{2}v(x),$ where $c$ is sufficiently   large positive number, is a strictly $(\widetilde{\gamma}_{2}+\varepsilon_{0})$-admissible subsolution to the problem \eqref{7.19},\eqref{7.2}. 
Then we can apply the result of Subsection 7.1 in the case $k=2$ and $B(x,z,p)=0$ to conclude
  the unique solvability of the problem \eqref{7.19},\eqref{7.2} in the class of strictly $\left(\widetilde{\gamma}_{2}+\frac{\varepsilon_{0}}{4}\right)$-admissible solutions. In this case the matrices $B(x, z, p)$ need not to be sufficiently small as in \eqref{7.17}, \eqref{7.18}, they satisfy only the conditions \eqref{7.15}-\eqref{7.22} and must not to be bounded on $\mathcal{D}.$

\end{document}